\documentclass{llncs}
\usepackage[hmargin=1.6in,vmargin=1.25in]{geometry}

\usepackage[dvips]{graphics}
\usepackage{epsfig}
\usepackage{url}
\usepackage{amssymb, amsmath}
\usepackage{float}
\usepackage{theorem}
\usepackage{xypic}
\xyoption{all}
\CompileMatrices

\newcommand{\QQ}{\mathbb{Q}}
\newcommand{\ZZ}{\mathbb{Z}}
\newcommand{\NN}{\mathbb{N}}
\newcommand{\RR}{\mathbb{R}}
\newcommand{\CC}{\mathbb{C}}
\newcommand{\KK}{\mathbb{K}}

\newcommand{\mm}{\mathfrak{m}}
\newcommand{\mL}{\mathcal{L}}
\newcommand{\mA}{\mathcal{A}}

\newcommand{\PP}{\mathbb{P}}
\newcommand{\Div}{\mathrm{Div}}

\newcommand{\FF}{\mathbb{F}}
\newcommand{\Fq}{\mathbb{F}_q}

\newcommand{\Qp}{\mathbb{Q}_p}

\newcommand{\Qq}{\mathbb{Q}_q}
\newcommand{\Zq}{\mathbb{Z}_q}

\newcommand{\Supp}{\mathrm{Supp}}

\newcommand{\of}{\overline{f}}

\newcommand{\oC}{\overline{C}}

\newcommand{\oA}{\overline{A}}

\newcommand{\bproof}{\noindent {\scshape Proof: }}
\newcommand{\eproof}{\ \mbox{} \hfill $\square$\mbox{}\newline}
{
\theoremstyle{plain}
\theoremheaderfont{\normalfont\scshape}
\theorembodyfont{\ttfamily}
\newtheorem{algo}{Algorithm}
}

\newcommand{\ceil}[1]{\left\lceil #1 \right\rceil}

\newcommand{\ord}{\mathrm{ord}}

\newcommand{\Vol}{\mathrm{Vol}}

\newcommand{\Spec}{\mathrm{Spec}\, }

\newfloat{algorithm}{ptbh}{algorithm.log}

\newcommand{\AB}{{\mathbb{A}}}
\newcommand{\Lin}{{\mathrm{Lin}}}
\newcommand{\TT}{{\mathbb{T}}}
\newcommand{\mS}{{\mathcal{S}}}

\newcommand{\nd}{nondegenerate with respect to its Newton polytope}

\begin{document}

\title{Computing Zeta Functions of Nondegenerate Curves}
\author{W.~Castryck\inst{1}\thanks{Research assistant of the Fund for Scientific Research - Flanders (FWO - Vlaanderen)},
J.~Denef\inst{1} and F.~Vercauteren\inst{2}\thanks{Postdoctoral fellow of the Fund for Scientific Research - Flanders (FWO - Vlaanderen)}}

\institute{
Department of Mathematics, \\
University of Leuven, \\
Celestijnenlaan 200B, \\
B-3001 Leuven-Heverlee, Belgium \\
\email{\{wouter.castryck, jan.denef\}@wis.kuleuven.be}
\and
Department of Electrical Engineering \\
University of Leuven \\
Kasteelpark Arenberg 10\\
B-3001 Leuven-Heverlee, Belgium \\
\email{frederik.vercauteren@esat.kuleuven.be}}

\maketitle

\begin{abstract}

In this paper we present a $p$-adic algorithm to compute the zeta function of a
nondegenerate curve over a finite field using Monsky-Washnitzer cohomology.  The paper
vastly generalizes previous work since in practice all known cases, e.g.\ hyperelliptic, superelliptic
and $C_{ab}$ curves, can be transformed to fit the nondegenerate case.
For curves with a fixed Newton polytope, the property of being nondegenerate is
generic, so that the algorithm works for almost all curves with given
Newton polytope. For a genus $g$ curve over $\FF_{p^n}$, the expected running time
is $\widetilde{O}(n^3 g^6 + n^2 g^{6.5})$, whereas the space complexity amounts to
$\widetilde{O}(n^3 g^4)$, assuming $p$ is fixed.

\end{abstract}

{\bf Keywords:} nondegenerate curves,
zeta function, Monsky-Washnitzer cohomology, Kedlaya's algorithm, Newton polytope, toric geometry,
effective Nullstellensatz


\section{Introduction}

An important research topic in computational number theory is the determination of
the number of rational points on an algebraic curve $\oC$ over a finite field $\FF_{p^n}$.
More generally, one is interested in the computation of its Hasse-Weil zeta function
\[ Z_{\oC}(t) = \text{exp}\left(\sum_{k=1}^\infty \# \oC(\FF_{p^{nk}}) \frac{t^k}{k} \right) \quad \in \QQ[[t]],\]
which turns out to be a rational function \cite{Dwork} (and hence a finite, computable object)
that contains a huge amount of arithmetic and geometric information about $\oC$.
For instance, if one wants to use a cryptosystem based on the discrete logarithm problem
on the Jacobian variety $\text{Jac}(\oC)$, one should be able to compute the cardinality
of its set of rational points, which is fully determined by $Z_{\oC}(t)$.
Efficient point counting algorithms can also
provide important heuristical (counter)evidence for several conjectures
concerning the asymptotic behavior of the number of points on algebraic curves (see for
instance \cite{Poonen,Katz,Poonen2}).

Mainly because of its applications in cryptography, a significant amount of work has been done in the field of
elliptic curve point counting. This roughly resulted in two types of algorithms. Schoof developed
a so-called $\ell$-adic algorithm \cite{SEA}, using torsion points to determine the number of points modulo small primes $\ell \neq p$.
This algorithm has polynomial running time in the input size $\sim n \log p$.
On the other hand, Satoh invented a $p$-adic method \cite{Satoh}, using the Serre-Tate canonical lift of the curve.
Unlike Schoof's algorithm, its running time is exponential in $\log p$.
For fixed (small) $p$ however, it is much faster, especially due to several improvements made in the past few years (see \cite{VThesis} for
an overview).

Generalizing the above techniques to curves of any genus is a
nontrivial task, since both methods make explicit use of the
geometry of elliptic curves. Another concern is that the resulting
algorithm should also have a good time complexity in the genus $g$
of the input curve, as its size should now be measured as $\sim gn
\log p$. So far, all attempts using the $\ell$-adic approach yield
impractical algorithms for $g > 2$ (see \cite{GaudrySchost} for a
treatment of the $g=2$ case, see also
\cite{HuangIerardi,HuangWong,Pila}), but the $p$-adic story is more
successful. In 2001, Kedlaya found a non-obvious way to `generalize'
Satoh's method to hyperelliptic curves of any genus\footnote{This
was over finite fields of odd characteristic. The characteristic 2
case was treated in \cite{DVHyp2}.} \cite{Kedlaya}, using a rigid
analytical lift instead of the canonical lift. The big technical
tool behind Kedlaya's algorithm is Monsky-Washnitzer cohomology (see
\cite{MW1,MW2,MW3} and the survey by van der Put \cite{vanderput}).

A particularly nice aspect of Kedlaya's method is that there are no obvious theoretical obstructions for
generalizations to larger classes of curves. This observation soon resulted in point counting algorithms
for superelliptic curves \cite{GaudryGurel2}
and $C_{ab}$ curves \cite{DVCab}. In the present paper, we vastly
generalize the previous by presenting an algorithm that determines the zeta function of
so-called \emph{nondegenerate} curves.
These are curves in $\left(\mathbb{A}^1_{\FF_{p^n}} \setminus \{0\}\right)^2$ that are defined by a Laurent polynomial
$\of \in \FF_{p^n}[x^{\pm 1},y^{\pm 1}]$ that is \emph{\nd}. We refer to Section~\ref{nondegcurves} for the definition
but mention here already that this condition is satisfied for \emph{generically chosen} Laurent polynomials with given
Newton polytope.

The main result can be formulated as follows.

\begin{theorem} \label{main}
There exists a deterministic algorithm to compute the zeta function of a genus $g$ nondegenerate
curve over $\FF_{p^n}$ that requires $\widetilde{O}(n^3\Psi_t)$ bit-operations and $\widetilde{O}(n^3\Psi_s)$ space
for $p$ fixed. Here, $\Psi_t$ and $\Psi_s$ are parameters that depend on the Newton polytope of the input curve only;
for `most common' Newton polytopes, $\Psi_t = \widetilde{O}(g^{6.5})$ and $\Psi_s = \widetilde{O}(g^4)$.
\end{theorem}

For explicit formulas for $\Psi_t$ and $\Psi_s$
we refer to Theorem~\ref{maintheorem} (Section~\ref{complsec}). Recall that the Soft-Oh notation $\widetilde{O}$
neglects factors that are logarithmic in the input size.
The notion `most common' is not intended to be made mathematically exact. It just means that
the Newton polytope should not be shaped too exotically. We refer to Section~\ref{complsec} for
more details.

It is worth remarking that Kedlaya's method is not the only $p$-adic
point counting technique that is being investigated for higher
genus. In 2002, Mestre adapted his so-called AGM method to ordinary
hyperelliptic curves of any genus over finite fields of
characteristic two \cite{Mestre}; it has been optimized by Lercier
and Lubicz \cite{LercierLubicz}, while Ritzenthaler extended it to
non-hyperelliptic curves of genus three \cite{Ritzenthaler}. These
algorithms have running time $\widetilde{O}(n^2)$ (for fixed $p$ and
$g$) but are exponential in the genus. Another interesting approach
is to combine Kedlaya's ideas with Dwork's deformation theory
\cite{Dworkdeform}. This was first proposed by Lauder \cite{Lauder1}
and has been studied in more detail by Lauder himself
\cite{Lauder2}, Gerkmann \cite{Gerkmann} and Hubrechts, who recently
obtained a memory efficient version of Kedlaya's original algorithm
\cite{Hubrechts}. Independently, Tsuzuki used similar ideas for
computing certain one-dimensional Kloosterman sums \cite{Tsuzuki}.

The remainder of this paper is organized as follows: Section~\ref{nondegcurves}
recalls the definition of nondegenerate curves, illustrates that a wealth of
information is contained in the Newton polytope and ends with a
new result on the effective Nullstellensatz problem.  Section~\ref{cohomsec} contains
a novel method to explicitly compute a basis of the first Monsky-Washnitzer cohomology
group  and Section~\ref{Frobeniuslift} describes
an algorithm to lift the Frobenius endomorphism.
An algorithm to compute modulo exact differential forms is
given in Section~\ref{redsec}.  Section~\ref{commodesec} discusses
the simplifications when the curve is commode and monic.
Finally, Section~\ref{complsec} contains the detailed algorithm and
complexity estimates and Section~\ref{conclsec} concludes the paper.\\

\noindent \textbf{Preliminaries.} Instead of giving a concise r\'esum\'e of the cohomology theory
of Monsky and Washnitzer, we immediately refer to the survey by van der Put \cite{vanderput} (or to the short overviews
given in e.g. \cite{Kedlaya} or \cite{DVCab}). The idea
behind the present algorithm is then simply to compute all terms in
the Lefschetz fixed point formula \cite[Formula \textbf{(1.2)}]{vanderput} (or
\cite[Theorem 1]{Kedlaya} or \cite[Theorem 1]{DVCab}) modulo
a certain $p$-adic precision.\\

\noindent \textbf{Notations and conventions.} Throughout this article, $x$ and $y$ are fixed formal variables.
For any integral domain $R$
and any subset $\mathcal{S} \subset \RR^2$, we denote by $R[\mathcal{S}]$ the ring generated by the monomials that
are supported in $\mathcal{S}$, i.e.
\[ R[x^iy^j \ | \ (i,j) \in \mathcal{S} \cap \ZZ^2] \, .\]
For instance, $R[\NN^2]$ is just the polynomial ring $R[x,y]$,
$R[\ZZ^2]$ is the Laurent polynomial ring $R[x^{\pm 1}, y^{\pm 1}]$, and so on.
If $R$ is a complete DVR with local parameter $t$, and if $R[\mathcal{S}]$ is a finitely generated $R$-algebra, we denote
its $t$-adic completion by $R \langle \mathcal{S} \rangle$ and its \emph{weak completion} by
$R\langle \mathcal{S} \rangle^\dagger$, see \cite{vanderput} (or
\cite{Kedlaya} or \cite{DVCab}) for the definition. Finally, if $\KK$ is a field, $\overline{\KK}$
denotes a fixed algebraic closure.

When dealing with cones or polytopes in $\RR^2$, we will often
implicitly assume that they are of full dimension,
that is: they are not contained in a line. However, this will always be clear from the context. If there is possible doubt,
the condition will be stated explicitly.


\section{Nondegenerate Curves}\label{nondegcurves}

Let $\KK$ be an arbitrary field and denote with
$\TT^2_\KK = \Spec \KK[\ZZ^2]$ the two-dimensional algebraic torus over $\KK$.
Consider
\[ f(x,y) = \sum_{(i,j) \in \mS} f_{i,j} x^i y^j \in \KK[\ZZ^2] \]
with $\mathcal{S}$ a finite subset of $\ZZ^2$ and $f_{i,j} \in \KK \setminus \{0\}$ for all $(i,j) \in \mS$.
The set $\mS$ is called the \emph{support} of $f$.  Denote by $\Gamma = \Gamma(f)$ the convex hull in $\RR^2$ of
the points $(i,j) \in \mS$, it is called the {\em Newton polytope} of $f$.  The boundary
of $\Gamma$ is denoted by $\partial \Gamma$.   The faces of $\Gamma$ can be subdivided
according to their dimension: \emph{vertices}, \emph{edges} and $\Gamma$
itself. Let
$\gamma$ be an edge between the integral points $(a,b)$ and $(c,d)$, then the
\emph{arithmetic length} $l(\gamma)$ is defined as $l(\gamma) = \gcd(a - c, b - d) \in \NN \setminus \{0\}$.
Note that the number of integral points on $\gamma$ is equal to $l(\gamma)+1$.

\begin{definition} \label{nondegenerate}
Let $f(x,y) = \sum_{(i,j) \in \mS} f_{i,j} x^i y^j \in \KK[\ZZ^2]$ be a Laurent polynomial with
Newton polytope $\Gamma$.  For each face $\gamma$ of $\Gamma$, define
$f_\gamma(x,y) = \sum_{(i,j) \in \gamma \cap \ZZ^2} f_{i,j}x^i y^j$. Then $f$ is called \nd \
if for all faces $\gamma$,
the system of equations
\[ f_\gamma = x \frac{\partial f_\gamma}{\partial x} = y \frac{\partial f_\gamma}{\partial y} = 0 \]
has no solutions in the torus $\TT^2_{\KK}$
(that is,
there are no solutions in $\left( \overline{\KK} \setminus \{0\} \right)^2$).
\end{definition}

Before recalling the geometric meaning of this notion, we prove that a sufficiently generic Laurent polynomial with
given Newton polytope will be nondegenerate. This is well-known in the characteristic 0 case.

\begin{lemma} \label{lemmageneric}
Let $\Gamma \subset \RR^2$ be the convex hull of a set of points in $\ZZ^2$. Consider the map
\[ \varphi : \ZZ^2 \rightarrow \mathbb{A}_{\KK}^2 : (i,j) \mapsto (i,j).\]
Then the dimension of the affine subspace of $\mathbb{A}_{\KK}^2$ spanned by $\varphi(\Gamma \cap \ZZ^2)$ equals
$\dim \Gamma$.
\end{lemma}

\bproof
This is not entirely trivial if $\KK$ is of characteristic $p \neq 0$. As the $\dim \Gamma = 0$ case is obvious,
we first suppose that $\dim \Gamma = 1$. Take points $q_1 \neq q_2 \in \Gamma \cap \ZZ^2$ and suppose that $\varphi(q_1) =
\varphi(q_2)$. Then we must have that $q_2 = q_1 + p^ev$ for some $e \in \NN_0$ and some
nonzero $v \in \ZZ^2$ that is not divisible by $p$.
Because $\Gamma$ is convex, it also
contains $q_1 + v$, and definitely $\varphi(q_1) \neq \varphi(q_1 + v)$.

Now suppose $\dim \Gamma = 2$. Take points $q_1, q_2 \in \Gamma \cap \ZZ^2$ such that $\varphi(q_1) \neq \varphi(q_2)$.
Take a $q_3 \in \Gamma \cap \ZZ^2$ that is not in the span of $q_1$ and $q_2$, but suppose $\varphi(q_3)$ \emph{is} in
the span of $\varphi(q_1)$ and $\varphi(q_2)$, say
\[ q_3 = q_1 + k(q_2 - q_1) + p^ev\]
for some $e \in \NN_0$ and some nonzero $v \in \ZZ^2$ that is not divisible by $p$ and
linearly independent of $q_2 - q_1$. Note that although this
expansion is far from unique, there is a natural upper bound for $e$, so that we may assume that
it is maximal. Indeed, if we write
$q_3 - q_1 = (a_1,a_2)$ and $q_2 - q_1 = (b_1,b_2)$, then it is not hard to see that $p^e | b_2a_1 - a_2b_1 \neq 0$.
As a consequence, $\varphi(v)$ and $\varphi(q_2 - q_1)$ are linearly independent, since otherwise this would contradict the maximality of $e$.

Next, we may suppose that $0 \leq k < p^e$
by repeatedly replacing $p^ev \leftarrow p^ev \pm p^e(q_2 - q_1)$ if necessary. We may even suppose that $k \neq 0$, since
otherwise we can proceed as in the $\dim \Gamma = 1$ case. Now define
\[ q = \frac{k-1}{p^e}q_1 + \frac{p^e-k}{p^e}q_2 + \frac{1}{p^e}q_3 = q_2 + v.\]
The first equality shows that $q \in \Gamma$, the second one shows that $q \in \ZZ^2$. Finally, $\varphi(q)$ is not
in the span of $\varphi(q_1)$ and $\varphi(q_2)$.
\eproof

\begin{proposition} \label{generic}
Let $\Gamma$ be a convex polytope in $\RR^2$ with integral vertex coordinates and write $\mathcal{S} = \Gamma \cap \ZZ^2$. Then the set of points
\[ (f_{i,j})_{(i,j) \in \mathcal{S}} \in \mathbb{A}_\KK^{\# \mathcal{S} } \]
for which $f = \sum f_{i,j} x^iy^j$ is \emph{not} \nd \ is contained in an algebraic set of codimension $\geq 1$.
Moreover, this algebraic set is defined over the prime subfield of $\KK$.
\end{proposition}

\bproof
Let $\gamma$ be a face of $\Gamma$. Suppose for now that it is two-dimensional.
Let $X_\gamma$ be the algebraic set in $\mathbb{A}_\KK^{\# \mathcal{S}}
\times \left(\mathbb{A}_{\KK} \setminus \{0\} \right)^2$
defined by the equations
\[ \sum_{(i,j) \in \gamma \cap \ZZ^2} f_{i,j} x^iy^j  =  0, \quad \sum_{(i,j) \in \gamma \cap \ZZ^2}
i f_{i,j} x^iy^j  =  0,
\quad \sum_{(i,j) \in \gamma \cap \ZZ^2} j f_{i,j} x^iy^j = 0.\]
It has codimension $3$. Indeed, for every $a,b \in \KK \setminus \{0\}$ the above equations define a linear codimension $3$ subspace of
$\mathbb{A}_\KK^{\# \mathcal{S}} \times \{ x=a,y=b \}$. Here we used that there is no
$(a,b,c) \in \KK^3 \setminus \{(0,0,0)\}$ such that $a + bi + cj = 0$ for all
$(i,j) \in \varphi(\gamma \cap \ZZ^2)$, where $\varphi$ is the map from the foregoing lemma.
Let $Y_\gamma$ be the projection of $X_\gamma$ on $\mathbb{A}_\KK^{\# \mathcal{S}}$. It has codimension at least 1 and consists
exactly of those $(f_{i,j})_{(i,j) \in \mathcal{S}}$ that correspond to a Laurent polynomial for which the
nondegenerateness condition with respect to $\gamma$ is not satisfied.

If $\gamma$ has dimension $< 2$, one can again construct such a $Y_\gamma$ using an appropriate change of variables
so that $f_\gamma$ becomes a univariate Laurent polynomial, or a constant.

Then the Zariski closure of $\cup_\gamma Y_\gamma$ is the requested algebraic set. Remark that $\cup_\gamma Y_\gamma$ may
contain points that correspond to Laurent polynomials that \emph{are} nondegenerate with respect to their Newton polytope:
this will be the case whenever they have a Newton polytope that lies strictly inside $\Gamma$.
\eproof

\begin{corollary}
Let $\Gamma$ be a convex polytope in $\RR^2$ with integral vertex coordinates and let $p$ be a prime number. Let $P_n$ be the probability that
a randomly chosen $\of \in \FF_{p^n}[\ZZ^2]$ with support inside $\Gamma$ is \nd. Then $P_n \rightarrow 1$
as $n \rightarrow \infty$.
\end{corollary}

Note that Proposition~\ref{generic} is false if
the condition of $\Gamma$ being convex is omitted: $\mathcal{S} = \{ (0,0), (p,0), (0,p) \}$
is an easy counterexample (where $p > 0$ is the field characteristic). Another important remark is
that Proposition~\ref{generic} cannot
be generalized to higher dimensions. For instance, \emph{any} trivariate polynomial having
\[ \Gamma =
\text{Conv}\{(0,0,0),(1,0,0),(0,1,0),(1,1,p)\}\] as its Newton polytope
(where $p>0$ is again the field characteristic) will have a singular point
in the three-dimensional algebraic torus.

Clearly, when $f$ is \nd,
then $f(x,y) = 0$ defines a non-singular curve on the torus $\TT^2_\KK$ (at least if $\dim \Gamma= 2$). But nondegenerateness is much stronger:
it implies that there exists a natural compactification $X_{\Gamma}$
of $\TT^2_\KK$ in which the closure of this curve is still non-singular.


\subsection{Toric Resolution} \label{TorRes}

The construction of $X_\Gamma$ is based
on the theory of toric varieties. We refer to \cite{Danilov} for the general theory. For the convenience of
the reader, we will explain the needed material in a self-contained way.

To any cone $\Delta \subset \RR^2$, i.e.\ the set of linear combinations with non-negative real
coefficients of a finite number of vectors in $\QQ^2$, we associate the \emph{affine toric surface}
$X_\Delta = \Spec \KK[\Delta]$.


Let $\Gamma$ be a polytope in $\RR^2$, then we can associate a \emph{toric surface} $X_\Gamma$
to $\Gamma$ in the following way: to each face $\gamma$, associate the
cone $\Delta(\gamma)$ generated by all vectors in
\[ \{x - p \ | \ x \in \Gamma, p \in \gamma \}  \, . \]
Let $U_\gamma$ be the affine toric surface $X_{\Delta(\gamma)}$.
If $\gamma \subset \tau$ with $\tau$ another face of $\Gamma$,
then $\Delta(\gamma) \subset \Delta(\tau)$ and $\KK[\Delta(\tau)]$ is obtained
from $\KK[\Delta(\gamma)]$ by adjoining the inverse of each monomial $x^iy^j \in \KK[\Delta(\gamma)]$ for which
$(i,j) \in \Lin(\tau)$.  Here $\Lin(\tau)$ is the linear subspace of
$\RR^2$ generated by the differences of vectors in $\tau$. Thus, $\Spec \KK[\Delta(\tau)]$ is obtained from
$\Spec \KK[\Delta(\gamma)]$ by cutting away some zero locus.
Otherwise said: $U_\tau$ is canonically embedded in $U_\gamma$ as a Zariski-open subvariety.
Note that $U_\Gamma = \TT^2_\KK$, so $\TT^2_\KK$ is canonically an open
subvariety of each variety $U_\gamma$.  The surface $X_\Gamma$ is then
covered by the affine toric surfaces $U_\gamma$ where $\gamma$ runs over all
vertices of $\Gamma$.  Two such surfaces $U_{\gamma_{1}}$ and $U_{\gamma_{2}}$
are glued together along their common open subvariety $U_\tau$ with
$\tau$ the smallest face of $\Gamma$ containing both $\gamma_1$ and $\gamma_2$.
The surface $X_\Gamma$ is complete and normal. Note
that the toric surface associated to (any multiple of) the standard $2$-simplex
is just the projective plane $\PP^2_\KK$.

To every face $\gamma$ we can associate the algebraic torus
$T_\gamma = \Spec \KK[\Lin(\gamma)]$.  Since $\Lin(\gamma) \subset \Delta(\gamma)$,
we obtain a surjective homomorphism from $\KK[\Delta(\gamma)]$
to $\KK[\Lin(\gamma)]$, by mapping the monomials $x^i y^j$ with
$(i,j) \in \Delta(\gamma) \setminus \Lin(\gamma)$ to zero and the other
monomials to themselves.  This canonically identifies $T_\gamma$ with
a closed subvariety of $U_\gamma$.  Note that $\dim T_\gamma = \dim \gamma$
and that $X_\Gamma$ is the disjoint union of the algebraic tori $T_\gamma$,
with $\gamma$ running over all faces of $\Gamma$.  Furthermore, the closure
of $T_\gamma$ in $X_\Gamma$ is the disjoint union of all the $T_\tau$
with $\tau$ a face of $\gamma$. Although $X_\Gamma$ may have singularities,
it is smooth outside the zero-dimensional locus associated
to the vertices of $\Gamma$.

Now, let $f(x,y) \in \KK[\ZZ^2]$ be a Laurent polynomial and let $\Gamma$ be its
Newton polytope. Let $V(f)$ denote the closure in $X_\Gamma$ of the locus of $f$
in the torus $\TT^2_\KK$, then $V(f)$ is called the {\em toric resolution} of
the affine curve defined by $f(x,y)=0$
on $\TT^2_\KK$.  Restricting $V(f)$ to $U_\gamma$, it is easy to verify
that $V(f) \cap T_\gamma$ equals the locus of $x^{-i} y^{-j} f_\gamma$ in $T_\gamma$,
where $(i,j) \in \gamma \cap \ZZ^2$.  A standard calculation then shows that
$V(f)$ intersects the torus $T_\gamma$ transversally and
does not contain $T_\tau$ for any vertex $\tau$ of $\gamma$ if $f_\gamma$, $\partial f_\gamma/ \partial x$
and $\partial f_\gamma / \partial y$ have no common zero in $\TT^2_\KK$.

In conclusion: the toric compactification $X_\Gamma$ of $\TT^2_\KK$ can be written
as the disjoint union
\begin{equation}
 X_\Gamma = \TT^2_\KK \cup T_1 \cup \cdots \cup T_r \cup P_1 \cup \cdots \cup P_r
\end{equation}
with $r$ the number of edges (and thus also the number of vertices) of $\Gamma$, $T_k$
the one-dimensional torus associated to the $k^\text{th}$ edge and
$P_k$ the zero-dimensional torus associated to
the $k^\text{th}$ vertex.  If $f$ is \nd \ $\Gamma$, then $V(f)$ is a complete nonsingular
curve on $X_\Gamma$ that does not contain $P_k$ for $k = 1, \ldots, r$ and
intersects the tori $T_k$ transversally for all $k$.


\subsection{Riemann-Roch and the Newton Polytope} \label{RRoch}

Most results in this subsection are easy consequences of known more general theorems \cite{Danilov}.
For the convenience of the reader we will give a self-contained exposition.
Throughout, assume that
$\KK$ is perfect. Let $f \in \KK[\ZZ^2]$ be \nd \ $\Gamma$ and let $C = V(f) \subset X_\Gamma$ be the toric
resolution of the curve defined by $f$ on $\TT^2_\KK$. Enumerate the vertices $p_1, \dots, p_r$
clockwise and let $t_k$ be the edge connecting $p_k$ with $p_{k+1}$ (where $p_{r+1} = p_1$).
Let $P_k \subset X_\Gamma$ be the zero-dimensional torus corresponding to $p_k$ and let $T_k \subset X_\Gamma$
be the one-dimensional torus corresponding to $t_k$.

For each $t_k$, denote with $e_k$ the vector
$(a_k,b_k) \in \ZZ^2$ with $\gcd(a_k, b_k) = 1$ which is perpendicular
to $t_k$ and points from $t_k$ towards the interior of $\Gamma$.
Define $N_k = p_k \cdot e_k$. Note that instead of $p_k$ we could have
taken any vertex on $t_k$, since the difference
is perpendicular to $e_k$.

Define the divisors $D_{C, \Gamma}$ and $W_C$ as
\[ D_{C, \Gamma} = - \sum_{k = 1, \ldots, r} N_k (T_k \cap C)
\quad \text{ and } \quad W_{C} =  \sum_{k = 1, \ldots, r} (T_k \cap C) \, . \]
The notation $D_{C, \Gamma}$ emphasizes that this divisor not
only depends on $C$, but also on $\Gamma$.
Indeed if we replace $f$ by $x^i y^j f$, then
$\Gamma$ is replaced by $\Gamma + (i,j)$, but $C$ remains the same.
Since from the context it will always be clear what $\Gamma$ is, we will mostly write $D_{C}$ instead
of $D_{C, \Gamma}$.

For any subset $A \subset \RR^2$, denote with $L_A$ the $\KK$-vector space generated
by $x^i y^j$ with $(i,j) \in A \cap \ZZ^2$.
If $D$ is a divisor on $C$ which is defined over $\KK$, then
$\mL(D)$ denotes the corresponding Riemann-Roch space
\[ \left\{ \left. f \in \KK(C) \setminus \{0\} \, \right| \, (f) + D \geq 0 \right\} \cup \{ 0 \} .\]
Note that $D_{C, \Gamma}$ and $W_C$ are defined over $\KK$. If $\KK \subset \KK'$ is a field extension,
we write
\[ \mathcal{L}_{\KK'}(D) = \left\{ \left. f \in \KK'(C) \setminus \{0\} \, \right| \, (f) + D \geq 0 \right\} \cup \{ 0 \}.\]
Note that $\mathcal{L}_{\overline{\KK}}(D)$
is generated by $\mathcal{L}(D)$ because $\KK$ is perfect.
In particular we have that $\dim_{\KK} \mathcal{L}(D)$ $= \dim_{\overline{\KK}} \mathcal{L}_{\overline{\KK}}(D)$.

We will often abuse notation
and write things as $L_A \subset \mL(D)$, though the latter is defined as a subspace of
the function field $\KK(C)$.

\begin{lemma} \label{fundamental}
Let $g$ be a Laurent polynomial with support in $m \Gamma$ for some $m \in \NN_0$. Let $P$ be
a point in $C \setminus \TT_\KK^2$ and denote with $t_k$ the edge of $\Gamma$ such that $P \in T_k$.
Then we have:
\begin{enumerate}
\item $\ord_P(g) \geq - \ord_P (m D_C)$;
\item If $g_{m t_k} = f_{t_k} = 0$ has no solutions in $\TT^2_\KK$, then equality holds. Conversely, if
equality holds for all $P \in T_k$, then $g_{m t_k} = f_{t_k} = 0$ has no solutions in $\TT^2_\KK$.
\end{enumerate}
\end{lemma}
\bproof
Let $p_k + \alpha$ be the integral point
on $t_k$ that is closest (but not equal) to $p_k$. Let $e_k = (a_k, b_k)$ be
as above, then $\alpha = (-b_k, a_k)$. Choose a vector $\beta = (c, d)$
such that
\[ \det \begin{pmatrix} -b_k & \ a_k \\ \ c & \ d \\ \end{pmatrix} = -1.\]
Note that the cone $\Delta(t_k)$ is generated by $\alpha, -\alpha, \beta$, so that
$U_{t_k} = \Spec \KK[x',x'^{-1},y'] \cong \mathbb{A}^2_\KK \setminus \mathbb{A}^1_\KK$. Here
\[ \begin{array}{ccc} x' & = & x^{-b_k}y^{a_k} \\ y' & = & x^cy^d \\ \end{array} \]
and $T_k$ corresponds to the locus of $y' = 0$ minus the origin.
Since $C$ intersects $T_k$ transversally, we have
that $y'$ is a local parameter for $C$ at $P$. Also note that $x'$ is a unit in the local ring at $P$.
The inverse transformation is given by
\begin{equation} \label{invtransform}
  \begin{array}{ccc} x & = & x'^{-d}y'^{a_k} \\ y & = & x'^cy'^{b_k} \\ \end{array}
\end{equation}
so that, using the notation $e'_k = (-d,c)$,
\begin{equation} \label{xiyj}
  x^iy^j = x'^{e_k' \cdot (i,j)}y'^{e_k \cdot (i,j)}.
\end{equation}
Now if $(i,j) \in m \Gamma$, then $e_k \cdot (i,j) \geq m (e_k \cdot p_k)$, with equality
iff $(i,j) \in m t_k$. Hence
\begin{equation} \label{gxy}
 g(x,y) = y'^{m(e_k \cdot p_k)} ( g_{m t_k} (x'^{-d}, x'^c) + y'( \cdots)).
\end{equation}
Since $e_k \cdot p_k = - \ord_P D_C$, the assertions follow. Indeed, $f_{t_k}(x'^{-d}, x'^c)$
vanishes at $P$, because (\ref{gxy}) also holds for $g$ replaced by $f$ and $m=1$.
\eproof

\begin{corollary} \label{inclusion} \mbox{} \\[-0.5cm]
\begin{enumerate}
\item For $i, j \in \ZZ$, the following holds:
\[ \Div_{C} (x^i y^j) = \sum_{k = 1, \ldots, r} (i,j) \cdot e_k (T_k \cap C) \, , \]
which implies that $L_{m\Gamma} \subset \mL(m D_C)$ for any $m \in \mathbb{N}_0$.
\item The arithmetic length $l(t_k)$ equals $\# ( T_k \cap C)$.
\end{enumerate}
\end{corollary}
\bproof
The first statement follows immediately from (\ref{xiyj}) and the
inequality on the line below it. The second
statement follows from the last assertion in the above proof,
namely that the points of $T_k \cap C$ correspond
to the zeros of $f_{t_k}(x'^{-d},x'^c)$. Now the latter can be written as
a power of $x'$ times a degree $l(t_k)$ polynomial in $x'$ with non-zero
constant term and without multiple roots.
\eproof

\begin{corollary} \label{canonicaldivisor}
Let $f_y$ denote the partial derivative of $f$ with respect to $y$, then
\[ \Div_C \left( \frac{dx}{x y f_y} \right) = D_{C} - W_C \, . \]
In particular, the differential form $dx / (xy f_y)$ has no poles, nor zeros
on $C \cap \TT^2_\KK$.
\end{corollary}
\bproof
First, let $P$ be a point of $C \setminus \TT^2_\KK$. We have to prove that
\[ \ord_P \frac{dx}{xyf_y} = \ord_P D_C - 1 \, . \]
With the notation as in Lemma \ref{fundamental}, we have that $f_{t_k}(x'^{-d}(P),x'^c(P))=0$, where
$k$ is such that $P \in T_k$.
Thus, because of the nondegenerateness of $f$:
\[ \left( x \frac{\partial f_{t_k}}{\partial x} \right) (x'^{-d}(P), x'^c(P)) \neq 0
  \quad \text{or} \quad \left( y \frac{\partial f_{t_k}}{\partial y} \right) (x'^{-d}(P), x'^c(P)) \neq 0 \, .\]
We may suppose that the second condition holds. Indeed, the first
case is treated analogously using that $dx / xyf_y = - dy/ xyf_x$.
Moreover, $\ord_P x$ is not a multiple of the characteristic $p$ of
$\KK$. Indeed if it was, then from formulas (\ref{invtransform}) and
the material above it, $a_k \equiv 0$ mod $p$ and $\alpha \equiv
(-b_k, 0)$ mod $p$ (if $p=0$, these congruences become exact
equalities). Hence $f_{t_k}$ has a special form: it equals a
monomial with exponent $p_k$ times a Laurent polynomial with all
exponents of $y$ divisible by $p$. This Laurent polynomial vanishes
on $(x'^{-d}(P),x'^c(P))$, because $x'$ is a unit at $P$. But this
contradicts the assumed second condition on $\frac{\partial
f_{t_k}}{\partial y}$.


Now apply Lemma \ref{fundamental} (and its proof) with $g$ replaced by $yf_y$ to find that
$\ord_P yf_y = - \ord_P(D_C)$. Since $\ord_P x$ is not divisible by $p$, we have that $\ord_P dx / x = -1$ and the result follows.

Next, take $P \in C \cap \TT^2_\KK$. Write $P = (p_x, p_y)$. Because of the
nondegenerateness we have that $\frac{\partial f}{\partial x}(P) \neq 0$ or
$\frac{\partial f}{\partial y}(P) \neq 0$. In particular, $dx / xyf_y = - dy/xyf_x$ can have no pole at $P$.
For the same reason, $x-p_x$ or $y-p_y$ must be local parameters at $P$ so that
for instance $dx/xyf_y = d(x-p_x)/xyf_y$ can have no zero at $P$.
\eproof

As a consequence, $D_C - W_C$ is
a canonical divisor. This observation allows us to give an elementary proof of a
well-known result. See \cite{Khovanskii} for much more general theorems on this matter.

\begin{corollary} \label{degDC} \mbox{} \\[-0.5cm]
\begin{enumerate}
\item $2g - 2 = \deg D_{C} - \deg W_{C}$, with $g$ the genus of $C$.
\item $g = \# ( (\Gamma \setminus \partial \Gamma) \cap \ZZ^2 )$, i.e.\
the genus of $C$ is the number of interior lattice points in $\Gamma$.
\end{enumerate}
\end{corollary}
\bproof
1. From the Riemann-Roch theorem it follows that the degree of a canonical
divisor is $2g-2$.

2. Because of Pick's theorem \cite{Pick}, which states that
\[ \Vol( \Gamma) = \# ( (\Gamma \setminus \partial \Gamma) \cap \ZZ^2 ) + \frac{ \# ( \partial \Gamma \cap \ZZ^2)}{2} - 1, \]
 it suffices to prove that $\deg D_C = 2 \Vol (\Gamma)$.
For every edge $t_k$, consider the triangle $\Delta_k$ defined by the two vertices of $t_k$ and the
origin. If the origin happens to be one of the vertices, this is just a line segment.
Then
\[ \Vol (\Gamma) = \sum_k - \text{sgn}(N_k) \text{Vol}(\Delta_k).\]
Now $\Delta_k$ is a triangle with base $l(t_k) \|e_k\|$ (the length of $t_k$)
and height $| p_k \cdot e_k |/ \|e_k\|$, so that its volume equals $l(t_k) |N_k| / 2$. The
result follows.
\eproof

We note that the inequality $g \leq \# ( (\Gamma \setminus \partial \Gamma) \cap \ZZ^2 )$
holds in any case, i.e. without the nondegenerateness condition. This is \emph{Baker's formula}: over $\CC$ it
was known already in 1893, 
a proof of the general case can be found in \cite{BeelenPellikaan}.

We conclude this subsection with the following theorem, which is an easy consequence of the fact
that $H^i(X_\Gamma, \mathcal{E}) = 0$ for any $i \geq 1$ and any invertible sheaf $\mathcal{E}$ on $X_\Gamma$
which is generated by its global sections (see \cite[Corollary~7.3 and Proposition~6.7]{Danilov}).
But for the convenience of the reader we will give a more elementary proof.

\begin{theorem} \label{RRSpace}
For any $m \in \mathbb{N}_0$, the Riemann-Roch space
$\mL(m D_C)$ is equal to $L_{m\Gamma}$.
\end{theorem}
\bproof
For our proof, the abuse of notation mentioned at the beginning of this subsection
is somewhat annoying. Therefore, we will temporarily introduce the notation $\mA_m$,
which denotes the image of $L_{m\Gamma}$ \emph{inside the function field $\KK(C)$}.
Note that the actual statement of the theorem should then be: $\mL(m D_C) = \mA_m$.

We already showed that $\mA_m \subset \mL(mD_C)$. Therefore, it suffices to
prove that the dimensions are equal. From Corollary \ref{degDC} and
the Riemann-Roch theorem, we have that $\dim \mL(mD_C) = m \deg D_C + 1 - g$.
Note that this is a polynomial of degree 1 in $m \geq 1$. Now consider the maps
\[ r_m : L_{(m-1)\Gamma} \rightarrow L_{m\Gamma} : w \mapsto wf, \quad m \geq 1 \, .\]
We claim that
$ \text{coker} \ r_m \cong \mA_m$. Indeed, we will show that the natural map
\[ \text{coker} \ r_m \rightarrow \mA_m \]
is injective. Let $v \in L_{m\Gamma}$ be such that $v=0$ in
the function field. Then there exists a unique Laurent polynomial $q$ such
that $v=fq$. Now for any $k=1, \dots, r$, we have that
\[ \ord_{T_k} v = \ord_{T_k} f  +  \ord_{T_k} q. \]
Here, $\ord_{T_k}$ is the valuation at $T_k$ in $X_\Gamma$ (which is nonsingular in codimension one).
From formula (\ref{gxy}) one deduces that $\ord_{T_k} v \geq m N_k$ (indeed, $y'$ is a local parameter
at $T_k$). Similarly, we have that $\ord_{T_k} f = N_k$. Therefore,
$\ord_{T_k} q \geq (m-1)N_k$. By a similar argument, now using (\ref{xiyj}), we conclude that
$q \in L_{(m-1)\Gamma}$, which proves the claim. Now by a well-known result by Ehrhart \cite{Ehrhart},
$\dim L_{m\Gamma}$ is a quadratic polynomial in $m$ with leading coefficient $\text{Vol}(\Gamma)$ for $m \geq 0$. As a consequence,
\[ \dim \mA_m = \dim \text{coker } r_m = \dim L_{m\Gamma} - \dim L_{(m-1)\Gamma} \]
is just like $\dim \mL(mD_C)$ a linear polynomial in $m$ for $m \geq 1$. Therefore it
suffices to prove equality for $m=1$ and $m \rightarrow \infty$.

The case $m=1$ follows from Corollary \ref{degDC}. Indeed,
\[ \dim \mL(D_C) = \deg D_C + 1 - g = 2g - 2 + \# (\partial \Gamma \cap \ZZ^2) + 1 - g = \# (\Gamma \cap \ZZ^2) - 1 \, ,\]
which is precisely $\dim \mA_1$.

For the case $m \rightarrow \infty$ it suffices to prove that
\[ \deg D_C = \lim_{m \rightarrow \infty} \frac{ \dim L_{m\Gamma} - \dim L_{(m-1)\Gamma} }{m} \, .\]
Since $\dim L_{m\Gamma} = \text{Vol}(\Gamma)m^2 + \dots$, the right hand side is $2 \text{Vol}(\Gamma)$
which is indeed
\[2 \# (\Gamma \setminus \partial \Gamma \cap \ZZ^2) + (\partial \Gamma \cap \ZZ^2) - 2
 = 2g - 2 + \deg W_C \]
according to Pick's theorem \cite{Pick}.
\eproof


\subsection{Effective Nullstellensatz}

In this subsection, we prove a new sparse effective Nullstellensatz.
Because this is interesting in its own right, things are treated somewhat more generally
than is needed for the rest of the paper. Let
$\KK$ be a field or a discrete valuation ring with maximal ideal $\mm$. Let $f \in \KK[\ZZ^n] := \KK[x_1^{\pm 1}, \dots, x_n^{\pm 1}]$ define
a smooth affine scheme over $\KK$. Then there exist Laurent polynomials
$\alpha, \beta_1, \dots, \beta_n \in \KK[\ZZ^n]$ for which
\[ 1 = \alpha f + \beta_1 x_1 \frac{\partial f}{\partial x_1} + \dots + \beta_n x_n \frac{\partial f}{\partial x_n}.\]
Though this is well known, we give the following inductive argument for the DVR case (using the field case),
for use in the proof of Lemma~\ref{DVRNullstellensatz}.
Let $t$ be a local parameter of $\KK$. Since $f$ is a smooth affine scheme over $\KK$, there exist Laurent polynomials
$\alpha, \beta_1, \dots, \beta_n \in \text{Frac}(\KK)[\ZZ^n]$ and $\tilde{\alpha}_1, \tilde{\beta}_1,
\dots, \tilde{\beta}_n \in \KK[\ZZ^n]$ such that
\begin{eqnarray}
\label{null1} 1 & = & \alpha f + \beta_1 x_1 \frac{\partial f}{\partial x_1} +
\dots + \beta_n x_n \frac{\partial f}{\partial x_n} \\
\label{null2} 1 & \equiv & \tilde{\alpha}f + \tilde{\beta}_1 x_1 \frac{\partial f}{\partial x_1} +
 \dots + \tilde{\beta}_n x_n \frac{\partial f}{\partial x_n} \mod t. \end{eqnarray}
Clearing denominators in (\ref{null1}) yields
\begin{equation} \label{null3}
  t^m = \alpha' f + \beta_1' x_1 \frac{\partial f}{\partial x_1} +
\dots + \beta_n' x_n \frac{\partial f}{\partial x_n}
\end{equation}
for some $m \in \mathbb{N}$ and $\alpha', \beta_1', \dots, \beta_n' \in \KK[\ZZ^n]$.
If $m=0$ we are done. If not, we can reduce $m$ inductively by rewriting
(\ref{null2}) as $1 + tQ = \tilde{\alpha}f + \tilde{\beta}_1 x_1 \frac{\partial f}{\partial x_1} +
\dots + \tilde{\beta}_n x_n \frac{\partial f}{\partial x_n}$, for some $Q \in \KK[\ZZ^n]$.
Now multiply this equation with $t^{m-1}$, multiply (\ref{null3}) with $Q$, and subtract. This concludes
the proof.

Next, for sake of being self-contained, we give the following definitions; they are straightforward generalizations
of the corresponding notions defined in a previous subsection.
\begin{definition}
The Newton polytope $\Gamma(h)$ of a polynomial $h \in \KK[\ZZ^n]$ is the polytope in $\RR^n$ obtained
by taking the convex hull of the support of $h$, i.e. the set of
exponent vectors in $\ZZ^n$ corresponding to monomials that have
a non-zero coefficient in $h$.
\end{definition}
For any $\sigma \subset \RR^n$, we denote by $h_\sigma$
the Laurent polynomial obtained from $h$ by setting all monomials whose exponent vectors lie outside of $\sigma$ equal to zero.
\begin{definition}
Suppose $\KK$ is a field and let $h \in \KK[\ZZ^n]$. Let $\Gamma$ be a convex polytope in $\RR^n$ with vertices in $\ZZ^n$ and
suppose $h$ has support inside $\Gamma$. If
\[ h_\gamma = x_1 \frac{\partial h_\gamma}{\partial x_1} = x_2 \frac{\partial h_\gamma}{\partial x_2}
= \dots = x_n \frac{\partial h_\gamma}{\partial x_n} = 0 \quad
\text{has no solutions in $\TT_\KK^n := \left( \overline{\KK}
\setminus \{0\} \right)^n$} \] for all faces $\gamma$ of $\Gamma$
(including $\Gamma$ itself), then it is said that $h$ is
nondegenerate with respect to $\Gamma$.
\end{definition}

Now, for reasons that will become clear in Section~\ref{Frobeniuslift}, we want the
Newton polytopes of $\alpha, \beta_1, \dots, \beta_n$ to be as small as possible. It will turn out that a natural
bound exists whenever $\of$, i.e.\ the reduction modulo $\mm$, is nondegenerate
with respect to $\Gamma(f)$, at least if the latter is $n$-dimensional and contains the origin.
The main theorem is the following sparse effective Nullstellensatz,
which seems new even when $\KK = \CC$. Our proof is inspired by an argument in \cite{Kushnirenko}, see also
\cite[Section 4]{Batyrevrepl}.





\begin{theorem} \label{Nullstellensatz}
Let $\KK$ be a field or a DVR, and denote
its maximal ideal by $\mm$. Let $\Gamma$ be a convex polytope in $\RR^n$
with vertices in $\ZZ^n$ and suppose that
$\dim \Gamma = n$. Let $f_0, f_1, \dots, f_n \in \KK[\ZZ^n]$ have
supports in $\Gamma$ and take a $g \in \KK[\ZZ^n]$ with support in $(n+1)\Gamma$.
Suppose that for every face $\gamma$
of $\Gamma$ the system $\of_{0 \gamma} =
\dots = \of_{n\gamma} = 0$ has no solutions in $\TT_{(\KK/\mm)}^n$.
Then there exist $h_0, \dots, h_n \in \KK[\ZZ^n]$ with
support in $n\Gamma$ such that $g = h_0f_0 + \dots + h_nf_n$.
\end{theorem}
\bproof
Write $k = \KK/\mm$. Let $S^k_\Gamma$ be the graded ring
consisting of all $k$-linear combinations of terms of
the form
\[ t^dx^e, \text{with $d \in \NN$ and $e \in d\Gamma \cap \ZZ^n$.}\]
The degree of such a term is by definition equal to $d$.
Similarly, let $S^\KK_\Gamma$ consist of the $\KK$-linear
combinations.

Let $\Delta$ be the cone in $\RR^{n+1}$ generated
by all vectors $(d,e)$ with $d \in \NN$ and $e \in d\Gamma$.
Clearly $S^k_\Gamma = k[\Delta]$. Because the systems
$\of_{0\gamma}= \dots= \of_{n\gamma} = 0$ have no common solution
in $\TT_k^n$, the locus in $\Spec(S^k_\Gamma)$ of
$(t\of_0, \dots, t\of_n)$ consists of only one point. This is easily
verified considering the restrictions of the locus of $tf_i$ to the tori
that partition $\Spec(S^k_\Gamma)$. Hence
\[ \frac{S^k_\Gamma}{(t\of_0, \dots, t\of_n)} \]
has Noetherian dimension zero. On the other hand $S^k_\Gamma$
is a Cohen-Macaulay ring by a well-known result of Hochster (see e.g.
\cite[Theorem 3.4]{Danilov}) that states that $k[C]$ is Cohen-Macaulay
for any cone $C$. So $t\of_0, \dots, t\of_n$ is a regular sequence.
This means that we have exact sequences
\[ 0 \rightarrow \left( \frac{S^k_\Gamma}{(t\of_0, \dots, t\of_i)} \right)_{d-1}
  \rightarrow \left( \frac{S^k_\Gamma}{(t\of_0, \dots, t\of_i)} \right)_d
  \rightarrow \left( \frac{S^k_\Gamma}{(t\of_0, \dots, t\of_{i+1})} \right)_d \rightarrow 0,\]
where the second arrow is multiplication by $t\of_{i+1}$ and
where $(\cdots)_d$ denotes the homogeneous part of degree $d$.
Thus
\[ \dim_k \left( \frac{S^k_\Gamma}{(t\of_0, \dots, t\of_{i+1})} \right)_d
   = \dim_k \left( \frac{S^k_\Gamma}{(t\of_0, \dots, t\of_i)} \right)_d
 - \dim_k \left( \frac{S^k_\Gamma}{(t\of_0, \dots, t\of_i)} \right)_{d-1}. \]
By a result of Ehrhart \cite{Ehrhart},
the number of lattice points in $d\Gamma$ (which is precisely $\dim_k (S_\Gamma)_d$) is a polynomial
function in $d$ for all $d\geq 0$. We obtain that
\[ \dim_k \left( \frac{S^k_\Gamma}{(t\of_0, \dots, t\of_n)} \right)_d\]
is a polynomial function in $d$ for all $d \geq n+1$. Since
the Noetherian dimension is zero, this polynomial must be zero as well.

In particular, we have that the $k$-linear map
\[ W_k : \bigoplus_{i=0}^n (S^k_\Gamma)_n \rightarrow (S^k_\Gamma)_{n+1} :
(t^n\overline{h}_0, \dots, t^n\overline{h}_n) \mapsto
t^{n+1}(\overline{h}_0\of_0 + \dots + \overline{h}_n\of_n)\]
is surjective. But then necessarily the corresponding $\KK$-map
\[ W_\KK : \bigoplus_{i=0}^n (S^\KK_\Gamma)_n \rightarrow
(S^\KK_\Gamma)_{n+1} : (t^nh_0, \dots, t^nh_n)
\mapsto t^{n+1}(h_0f_0 + \dots + h_nf_n)\]
is surjective. Indeed, let $M$ be the matrix of $W_\KK$. Then its
reduction modulo $\mm$ is the matrix of $W_k$, so it has a minor of maximal
dimension with non-zero determinant. But this means that $M$ itself
has a minor of maximal dimension whose determinant is a unit in $\KK$.
\eproof

The following corollary will be essential in devising a sharp bound for
the rate of overconvergence of a lift of the Frobenius endomorphism in Section~\ref{Frobeniuslift}. It
will also allow us to translate the action of Frobenius on the first Monsky-Washnitzer cohomology
space (consisting of \emph{differential forms})
to a space of \emph{functions} that will be introduced in the next section. The exact way in which
this is done is described in \verb"STEP I" and \verb"STEP V" of the algorithm (Section~\ref{complsec}).

\begin{corollary} \label{nullbound}
Let $\KK$ be a field or a DVR and denote by $k$ its residue field.
Let $f \in \KK[\ZZ^n]$ and suppose that $f$ and its reduction $\of
\in k[\ZZ^n]$ have the same Newton polytope $\Gamma$, which is
supposed to be $n$-dimensional and to contain the origin. If $\of$
is nondegenerate with respect to its Newton polytope, then there
exist $\alpha, \beta_1, \dots, \beta_n \in \KK[\ZZ^n]$ such that
\[ 1 = \alpha f + \beta_1 x_1 \frac{\partial f}{\partial x_1} + \dots + \beta_n x_n \frac{\partial f}{\partial x_n} \]
with $\Gamma(\alpha), \Gamma(\beta_1), \dots, \Gamma(\beta_n) \subset n \Gamma(f)$.
\end{corollary}
\bproof
Apply Theorem~\ref{Nullstellensatz} to
$f, x_1 \frac{\partial f}{\partial x_1}, \dots, x_n \frac{\partial f}{\partial x_n}$.
\eproof

A second corollary to Theorem~\ref{Nullstellensatz} is that an arbitrary lift
of a nondegenerate Laurent polynomial with the same Newton polytope is again nondegenerate.

\begin{corollary} \label{liftnondegenerate}
Let $\KK$ be a DVR with residue field $k$ and let $f \in \KK[\ZZ^n]$.
Suppose $f$ and its reduction $\of$ have the same Newton polytope. If
$\of$ is nondegenerate with respect to its Newton polytope, then so is $f$
(when considered over the fraction field of $\KK$).
\end{corollary}
\bproof
Let $\Gamma = \Gamma(f) = \Gamma(\of)$. Let $\gamma$ be any face of $\Gamma$.
If $\of$ is nondegenerate with respect to $\Gamma$, then so is $\of_\gamma$
with respect to $\gamma$. Using an appropriate change of variables, we can
apply Theorem \ref{Nullstellensatz} to find a Laurent monomial $x_1^{r_1} \dots x_n^{r_n}$ and Laurent
polynomials $g_0, g_1, \dots, g_n \in \KK[\ZZ^n]$ such that
\[ x_1^{r_1} \dots x_n^{r_n} = g_0 f_\gamma + g_1 \frac{\partial f_\gamma}{\partial x_1}
+ \dots + g_n \frac{\partial f_\gamma}{\partial x_n}.\]
In particular, $f_\gamma = x_1 \frac{\partial f_\gamma}{\partial x_1} =
\dots = x_n \frac{\partial f_\gamma}{\partial x_n} = 0$ can have no
solutions in $\TT_\text{Frac($\KK$)}^n$.
\eproof

We conclude this subsection with a discussion on what can happen if our Laurent polynomial is \emph{not} \nd.
In that case much worse bounds than the one given in Corollary~\ref{nullbound} need to be used.
We restrict our examples to the bivariate polynomial case, i.e. $f \in \KK[x,y]$.
If $\KK$ is a field, a quadratic upper bound follows from general effective Nullstellensatz
theorems, such as~\cite[Theorem 1.5]{Kollar}. Example~\ref{Nullstellenexample1} shows
that this bound is asymptotically sharp. If $\KK$ is a DVR, it is even impossible to give bounds
in terms of $\Gamma(f)$, which is shown in Example~\ref{Nullstellenexample2}.

\begin{example} \label{Nullstellenexample1}
Let $\KK$ be a field of finite characteristic $p$ and let $d \geq 2p$
be a multiple of $p$. Consider the degree $d+1$ polynomial
\[ f  =  y^{d+1} + x^{d-p}y^p + 1\]
(its definition is
inspired by \cite[Example 2.3]{Kollar}).
It defines an irreducible, nonsingular curve in $\mathbb{A}^2_{\KK}$,
so take polynomials $\alpha, \beta, \gamma \in \KK[x,y]$ such that
\[ 1 = \alpha f + \beta \frac{\partial f}{\partial x} + \gamma \frac{\partial f}{\partial y}.\]
Let $\lambda$ be the maximum of the degrees of $\alpha, \beta, \gamma$.
Homogenizing the above equation with respect to a new variable $z$ yields that
$z^{\lambda + d + 1} \in I$. Here, $I \subset \KK[x, y, z]$ is the ideal generated
by the homogenizations of $f$, $\frac{\partial f}{\partial x}$, $\frac{\partial f}{\partial y}$, i.e.
\[ I = (x^{d-p}y^pz + z^{d+1}, y^d ).\]
Now consider its image $\widetilde{I}$ under the map
\[ \KK[x, y, z] \rightarrow \KK[y, z] : h(x,y,z) \mapsto h(1, y, z). \]
Then
\[ \widetilde{I} = (y^pz + z^{d+1}, y^d).\]
It is easy to verify that no power of $z$ less than $d^2/p$ can be contained
in $\widetilde{I}$, and a fortiori the same holds for $I$. Therefore,
\[ \lambda \geq \frac{d^2}{p} - d - 1, \]
which is $O((\deg f)^2)$.
\end{example}

\begin{example} \label{Nullstellenexample2}
Let $\KK$ be an arbitrary DVR with local parameter $t$. Consider
$f = y - txy + (t^m + t^2)x^2 - 1 \in \KK[x,y]$ for some big
natural number $m$ (its definition is inspired by an example in
\cite{Asschenbrenner}). Since the system of equations
\[ \left\{ \begin{array}{cclcc} f & = & y - txy + (t^m + t^2)x^2 - 1 & = & 0 \\
                                \frac{\partial f}{\partial x} &=& -ty + 2(t^m + t^2)x          & = & 0 \\
                                \frac{\partial f}{\partial y} &=& 1 - tx                       & = & 0 \\ \end{array} \right. \]
has no solutions, neither over $\overline{\KK}$, nor over $\overline{\KK/(t)}$,
there exist polynomials $\alpha, \beta, \gamma \in \KK[x,y]$ that satisfy
\[ 1 = \alpha f + \beta \frac{\partial f}{\partial x} + \gamma \frac{\partial f}{\partial y}.\]
Putting $y = 1 + tx$ and reducing modulo $t^m$ gives the following identity in $\KK/(t^m)[x]$:
\[ 1 = \overline{\alpha} ( (1 - tx)(1 + tx) + t^2x^2 - 1) + \overline{\beta} (-t(1+tx) + 2t^2x) + \overline{\gamma} (1 - tx)\]
\[ = - t \overline{\beta} (1 - tx) + \overline{\gamma} (1 - tx)  = (\overline{\gamma} - t \overline{\beta})(1 - tx). \]
Since the inverse of $1-tx$ in $\KK/(t^m)[x]$ is $1 + tx + t^2x^2 + \dots + t^{m-1}x^{m-1}$, we conclude that
\[ \max \{ \deg \beta, \deg \gamma \} \geq \deg (\gamma - t\beta) \geq m - 1.\]
\end{example}

In fact, the above example shows that even the valuations of the coefficients of $f$ do not suffice to give a
Nullstellensatz bound. The best we can do are results of the following type.

\begin{lemma} \label{DVRNullstellensatz}
Let $\KK$ be a DVR with local parameter $t$. For every $d \in \mathbb{N} \setminus \{0\}$ there exists a non-zero polynomial
$g_d \subset \KK[c_{ij}]_{i,j \in \NN, i+j \leq d}$ of degree $\leq 3(d^2 + 1)(d^2 + 2)/2$,
for which the following holds. If
\[ f = \sum_{i + j \leq d} C_{ij} x^iy^j \quad \in \KK[x,y] \]
defines a smooth affine $\KK$-scheme, then there are polynomials
$\alpha, \beta, \gamma \in \KK[x,y]$ such that $1 = \alpha f + \beta \frac{\partial f}{\partial x} + \gamma \frac{\partial f}{\partial y}$ with
\[ \deg \alpha, \deg \beta, \deg \gamma \leq d^2 (1 + \ord_t g_d(C_{ij})). \]
\end{lemma}

\bproof
Given such an $f$,
we know from~\cite[Theorem 1.5]{Kollar} that there exist
$\alpha', \beta', \gamma' \in \text{Frac}(\KK)[x,y]$ of degree $\leq d^2$ for which
$1 = \alpha' f + \beta' \frac{\partial f}{\partial x} + \gamma' \frac{\partial f}{\partial y}$.
In other words, the formula
\[ 1 = \left(\sum_{i + j \leq d^2} \alpha'_{ij}x^iy^j\right) f +
\left(\sum_{i + j \leq d^2} \beta'_{ij}x^iy^j\right) \frac{\partial f}{\partial x} +
\left(\sum_{i + j \leq d^2} \gamma'_{ij}x^iy^j\right) \frac{\partial f}{\partial y} \]
gives rise to a system $\mathcal{S}_f$ of linear equations in $n = 3(d^2 + 1)(d^2 + 2)/2$ unknowns
$\alpha'_{ij}, \beta'_{ij}, \gamma'_{ij}$ that is solvable over $\text{Frac}(\KK)$. Let
\[ r := \max_f \text{rank}(\mathcal{S}_f) \leq n, \]
and let $f_0 = \sum_{i+j \leq d} C_{0,ij}x^iy^j$ be a polynomial for which this rank is actually obtained.
Then $\mathcal{S}_{f_0}$ has a non-zero $(r \times r)$-minor, which is a degree $r$
polynomial expression in the $C_{0,ij}$. Let $g_d(c_{ij}) \in \KK[c_{ij}]$ be the corresponding polynomial.

Now, using Cramer's rule, we can find a solution to $S_{f_0}$
such that the valuations of the denominators appearing in this solution are bounded by $\ord_t g_d(C_{0,ij})$.
In fact, this statement holds in general: for \emph{any} $f = \sum_{i+j \leq d} C_{ij}x^iy^j$, we can find a solution to $S_f$
whose denominators are bounded by $\ord_t g_d(C_{ij})$. Indeed, either $g_d(C_{ij})$ equals zero, or it is a minor of maximal dimension
of $S_f$. Now using the induction procedure mentioned at the beginning of this subsection, and again
using~\cite[Theorem 1.5]{Kollar} (but now over the residue field), we get the desired result.
\eproof


\section{Cohomology of Nondegenerate Curves}\label{cohomsec}

Let $\Fq$ be a finite field with $q = p^n$, $p$ prime.
Given a Laurent polynomial $\of \in \Fq[\ZZ^2]$
that is \nd \ $\Gamma$,
let $\oC$ denote the nonsingular curve $V(\of)$, i.e.\ the closure
in $\overline{X}_\Gamma$ of the zero locus of $\of$ in the torus $\TT^2_{\Fq}$.
Here $\overline{X}_\Gamma$ is the toric $\Fq$-surface
associated to $\Gamma$. Suppose
that $\overline{C}$ has genus $g \geq 1$.
Then without loss of generality we may assume that
\begin{enumerate}
  \item there are $d_t > d_b \in \ZZ$ such that
        $\Gamma$ has a unique top vertex (with $y$-coordinate $d_t$) and a unique bottom vertex
                 (with $y$-coordinate $d_b$)
  \item the origin is an
interior point of $\Gamma$.
\end{enumerate}
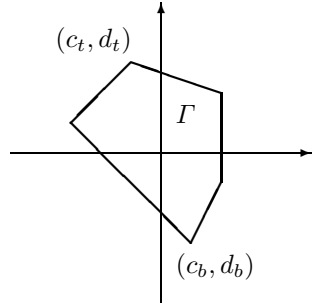
\begin{center}
\begin{picture}(20,20)

  \thinlines

  \put(0,10){\vector(1,0){20}}
  \put(10,0){\vector(0,1){20}}

  \put(3,17){\text{$(c_t,d_t)$}}
  \put(11,2){\text{$(c_b,d_b)$}}
  \put(11,12){\text{$\Gamma$}}

  \thicklines

  \put(8,16){\line(-1,-1){4}}
  \put(4,12){\line(1,-1){8}}
  \put(12,4){\line(1,2){2}}
  \put(14,8){\line(0,1){6}}
  \put(14,14){\line(-3,1){6}}
\end{picture}
\end{center}
The two important consequences of this setting are that
\begin{enumerate}
  \item the set $S := \{ x^ky^l \ | \ k,l \in \ZZ, \ d_b \leq l < d_t \ \}$ is an $\FF_q$-basis for $\frac{\FF_q[\ZZ^2]}{(\of)}$;
  \item \emph{every} $h \in \FF_q[\ZZ^2]$ has support in $m\Gamma$ for some big enough $m \in \mathbb{N}$.
\end{enumerate}
To see why we may assume that $\Gamma$ is of the above shape, take vertices $v_t$ and $v_b$
such that $\|v_t - v_b\|$ is maximal.
Take $\alpha \in \ZZ^2$ with coprime coordinates such that it is perpendicular to $v_t - v_b$.
Let $\beta \in \ZZ^2$ be such that $\alpha = (\alpha_1, \alpha_2), \beta = (\beta_1, \beta_2)$ generate $\ZZ^2$ over $\ZZ$.
There exist $(\gamma_1, \gamma_2), (\delta_1, \delta_2) \in \ZZ^2$ such that
\[ \left( \begin{array}{cc} \alpha_1 & \alpha_2 \\ \beta_1 & \beta_2 \\ \end{array} \right)
   \left( \begin{array}{cc} \gamma_1 & \gamma_2 \\ \delta_1 & \delta_2 \\ \end{array} \right) = \mathbb{I}. \]
Then the isomorphism $x \leftarrow x^{\gamma_1}y^{\gamma_2}$, $y \leftarrow x^{\delta_1}y^{\delta_2}$ yields
unique top and bottom vertices of the Newton polytope. Moreover, it preserves nondegenerateness.
For the second condition, remember that the number of interior points
in $\Gamma$ equals $g \geq 1$. Let $(a,b)$ be
an interior point, then $x^{-a} y^{-b} f$ is \nd \ $\Gamma - (a,b)$
and clearly $(0,0)$ is an interior point of $\Gamma - (a,b)$.

%
%
%
%
%
%
%
%
%
%

Let $\Qq$ denote the unramified extension of degree $n$ of $\Qp$
with valuation ring $\Zq$ and residue field $\Zq/ (p \Zq) \cong \Fq$.
Take an arbitrary lift $f \in \Zq[\ZZ^2]$ of $\of$ such that
$\Gamma(f) = \Gamma(\of) = \Gamma$. Note that we have properties similar to the ones mentioned above:
\begin{enumerate}
  \item the set $S := \{ x^ky^l \ | \ k,l \in \ZZ, \ d_b \leq l < d_t \ \}$ is a $\ZZ_q$-module basis for $\frac{\ZZ_q[\ZZ^2]}{(\of)}$;
  \item \emph{every} $h \in \ZZ_q[\ZZ^2]$ has support in $m\Gamma$ for some big enough $m \in \mathbb{N}$.
\end{enumerate}
Due to
Corollary \ref{liftnondegenerate}, $f$ is also \nd \ (when considered over $\QQ_q$). As a consequence, if $C$ denotes
the nonsingular curve obtained by taking the closure in $X_\Gamma$ (the toric $\QQ_q$-surface
associated to $\Gamma$) of the zero
locus of $f$ in $\TT^2_{\Qq}$, we have that $g(C) = g(\oC)$. Also,
for each edge $\gamma$ of $\Gamma$ we have
\[ \# (C \cap T_\gamma) = \# (\oC \cap \overline{T}_\gamma  ) \, , \]
where $T_\gamma$ (resp.\ $\overline{T}_\gamma$) is the algebraic
torus associated to $\gamma$ over $\Qq$ (resp.\ over $\Fq$) and
the intersections are transversal. These and other observations
reveal a deep geometric correspondence between $C$ and $\oC$ that
is directed by the Newton polytope. This is the main reason why we
work with nondegenerate curves. Indeed, as in Kedlaya's algorithm
it enables us to compute in the algebraic de Rham cohomology of
$C$. Namely, we have the following very general theorem
\cite[Theorem~1]{Kedlayaants}.

\begin{theorem}\label{comparisontheo}
Let $Y$ be a smooth proper $\Zq$-scheme, let $Z \subset Y$ be
a relative normal crossings divisor and let $X = Y \setminus Z$.
If $X$ is affine, then for any $i \in \mathbb{N}$ there exists a canonical isomorphism
\begin{equation} \label{isomorphism}
 H^i_{DR} (X \otimes_{\Zq} \Qq) \rightarrow H^i_{MW} (X \otimes_{\Zq} \Fq / \Qq) \, ,
\end{equation}
where $H_{DR}^i$ denotes algebraic de Rham cohomology and
$H_{MW}^i$ denotes Monsky-Washnitzer cohomology.
\end{theorem}

The above theorem applies in our situation with $X = \Spec \frac{\Zq[\ZZ^2]}{(f)}$ and $Y$
its closure in the toric scheme associated to $\Gamma$ (this is constructed exactly as in Section~\ref{TorRes},
with $\KK$ replaced by $\ZZ_q$). This is a smooth proper scheme and by the above observations
$Z = Y \setminus X$ is indeed a relative normal crossings divisor.

An alternative proof of Theorem~\ref{comparisontheo} (in the case
of a nondegenerate curve) follows from the material in Section~\ref{redsec}.
There we will implicitly prove that the canonical map (\ref{isomorphism}) with $i=1$ is
surjective. Since both $Z \otimes_{\Zq} \Fq$ and
$Z \otimes_{\Zq} \Qq$ contain $\# (\partial \Gamma \cap \ZZ^2)$ points,
we have that the dimensions are equal, which shows that the map is an isomorphism:
\[ \dim H^1_{DR} (X \otimes_{\Zq} \Qq) =
   \dim H^1_{MW} (X \otimes_{\Zq} \Fq / \Qq) =
   2g + \# (\partial \Gamma \cap \ZZ^2) - 1 = 2 \Vol(\Gamma) + 1 \, . \]
The last equality follows from Pick's theorem~\cite{Pick} and Corollary~\ref{degDC}.
This relationship between $\Vol(\Gamma)$ and the Betti numbers of $X \otimes_{\Zq} \Qq$ was
already noticed in a much more general setting by Khovanskii~\cite{Khovanskii}.

As a first step towards constructing a basis for $H^1_{DR}(X \otimes_{\Zq} \Qq)$,
we prove the following theorem. Let $A = \frac{\ZZ_q[\ZZ^2]}{(f)}$, denote by
$D^1(A)$ the universal $\ZZ_q$-module of differentials of $A$, and by $d: A \rightarrow D^1(A)$
the corresponding exterior derivation. Thus $H^1_{DR}(C \cap \TT^2_{\QQ_q}) = \frac{D^1(A)}{d(A)} \otimes_{\ZZ_q} \QQ_q$.

\begin{theorem} \label{reduction}
Every element of $D^1(A) \otimes_{\ZZ_q} \QQ_q$ is equivalent modulo
$d(A) \otimes_{\ZZ_q} \QQ_q$ with a differential form $\omega$ with divisor
\[ \Div_C(\omega) \geq - D_C - W_C \, . \]
\end{theorem}
\bproof
First, suppose that all places $P_k \in C \setminus \TT^2_{\QQ_q}$ are $\QQ_q$-rational.
Write $D_C = \sum_{k=1}^r a_k P_k$. Since the origin is an interior point of $\Gamma$, all $a_k >0$.
Note that $D_C + W_C = \sum_{k=1}^r (a_k + 1)P_k$ and that $\deg D_C = \sum_{k=1}^r a_k > 2g - 2$
due to Corollary~\ref{degDC}. Now let $\omega' \in D^1(A) \otimes_{\ZZ_q} \QQ_q$ have a pole of
order $b_k + 1 > a_k + 1$ at some place $P_k$. Because of the Riemann-Roch theorem,
we can find a function $h \in \mathcal{L}(a_1P_1 + \dots + b_kP_k + \dots + a_rP_r) \setminus
\mathcal{L}(a_1P_1 + \dots + (b_k - 1)P_k + \dots + a_rP_r)$. Then adding to $\omega'$
a suitable multiple of $dh$ will reduce the pole order at $P_k$. Continuing in this
way eventually proves the theorem in case all $P_k$ are $\QQ_q$-rational.

The general case follows easily from the above, using that $D_C$ and $W_C$ are defined over $\QQ_q$.
\eproof

In the above proof we reduced the pole orders one place at a time for simplicity.
However, in the algorithm the reduction will proceed more simultaneously by moving from
`level' $mD_C + W_C$ to `level' $(m-1)D_C + W_C$. Indeed, because of Theorem \ref{RRSpace}
we have a good understanding of what the spaces $\mL(mD_C)$ look like. For the same reason, it
is more natural to work with \emph{functions} instead of \emph{differential forms}.

Consider the following map
\[ \Lambda : A \otimes_{\ZZ_q} \QQ_q \rightarrow D^1(A) \otimes_{\ZZ_q} \QQ_q : h \mapsto h \frac{dx}{xyf_y} \, . \]
Let $\alpha, \beta \in A$ be such that $1 = \alpha f_x + \beta f_y$. Then $dx/f_y = \beta dx - \alpha dy$,
which shows that $\Lambda$ is well-defined. Moreover, for any
$g_1, g_2 \in A$ we have that $\Lambda(xy(f_yg_1 - f_xg_2)) = g_1dx + g_2dy$, so
that it is in fact a bijection.
By Corollary~\ref{canonicaldivisor} we have
$\Div_C (\Lambda(h)) = \Div(h) + D_C - W_C$.
Applying Theorem~\ref{reduction} shows that each differential form in
$D^1(A) \otimes_{\ZZ_q} \QQ_q$ is equivalent modulo exact differential forms
with an $\omega$ with divisor $\Div_C(\omega) \geq -D_C - W_C$.  Take
 $h \in A \otimes_{\ZZ_q} \QQ_q$ such that $\Lambda(h) = \omega$, then
\[ \Div_C(\Lambda(h)) = \Div_C(h) + D_C - W_C \geq - D_C - W_C \Leftrightarrow h \in \mL(2 D_C) \, . \]
This shows that $\Lambda( \mL(2 D_C))$ generates $H^1_{DR}(C \cap \TT^2_{\Qq})$.

To find an actual basis for $H^1_{DR}(C \cap \TT^2_{\Qq})$, we define an operator $D$ on $A$
such that the image under $\Lambda$ is an exact differential, i.e.\
$d h = \Lambda( D h)$.  The definition of $D$ follows easily from the following:
\[ d h = h_x dx + h_y dy = x y (f_y h_x  - f_x h_y) \frac{dx}{x y f_y}  \, , \]
and thus
\[ D(h) = x y \left( f_y \frac{\partial}{\partial x}  - f_x \frac{\partial}{\partial y} \right) (h)  \, . \]

By Theorem~\ref{RRSpace} we have the following corollary. A related description,
for nondegenerate hypersurfaces of any dimension, is contained in \cite[Corollary 6.10 and Theorem 7.13]{Batyrevrepl}.
\begin{corollary}
If the origin is an interior point of $\Gamma$, then $\Lambda$ induces an isomorphism of
$\Qq$-vector spaces:
\[ \frac{L_{2 \Gamma}}{f L_{\Gamma} + D ( L_\Gamma) } \cong  H^1_{DR}(C \cap \TT^2_{\Qq}) \, . \]
\end{corollary}
Note that the proof of the above corollary does not provide an explicit bound on the
denominators introduced during the reduction, which is required to
determine the $p$-adic precision up to which one has to compute.
In Section~\ref{redsec} we describe a simple
reduction algorithm and at the same time prove tight bounds on
the loss of precision.


\section{Lifting Frobenius Endomorphism} \label{Frobeniuslift}

We begin this section by introducing the following notation. Write
\[ \oA = \frac{\FF_q[\ZZ^2]}{(\of)} \, , \quad A = \frac{\ZZ_q[\ZZ^2]}{(f)} \, , \quad
A^\infty = \frac{\ZZ_q \langle \ZZ^2 \rangle}{(f)} \, , \quad A^\dagger = \frac{\ZZ_q \langle \ZZ^2 \rangle^\dagger}{(f)} \, ,\]
and define the $p^\text{th}$ and $q^\text{th}$ power Frobenius endomorphisms
\[ \overline{\mathcal{F}}_p : \oA \rightarrow \oA : a \mapsto a^p \, , \quad \overline{\mathcal{F}}_q : \oA \rightarrow \oA : a \mapsto a^q \,.\]
A main task in developing a point counting algorithm using Monsky-Washnitzer cohomology is the computation of a $\ZZ_p$-algebra endomorphism
$\mathcal{F}_p : A^\dagger \rightarrow A^\dagger$
that \emph{lifts} $\overline{\mathcal{F}}_p$ in the sense that $\overline{\mathcal{F}}_p \circ \pi = \pi \circ \mathcal{F}_p$, where
$\pi$ is reduction modulo $p$. Then $\mathcal{F}_q := \mathcal{F}_p \circ \dots \circ \mathcal{F}_p$ is a
$\ZZ_q$-algebra morphism that lifts $\overline{\mathcal{F}}_q$. Note that decomposing $\mathcal{F}_q$ into
$n$ copies of $\mathcal{F}_p$ ($p$ small) dramatically improves the running time of the algorithm: this is
the main reason why $p$-adic point counting algorithms are especially well-suited for small values of $p$.

First, we consider the following Hensel-like lemma.
In a paper subsequent to this one,
Kedlaya gives a related result, with a more elegant proof~\cite{Kedlayadraft}.

\begin{lemma} \label{hensel}
Let $\Gamma$ be a convex polygon in $\RR^2$ with vertices in $\ZZ^2$.
Take $a,b \in \mathbb{N}$ (not both zero) and let $H(Z) = \sum h_k Z^k \in \ZZ_q[\ZZ^2][Z]$ satisfy
\begin{enumerate}
  \item  $\Gamma(h_k) \subset (ak + b)\Gamma$ for all $k \in \NN$;
  \item $h_0 \equiv 0 \mod p$;
  \item $h_1 \equiv 1 \mod p$.
\end{enumerate}
Then there exists a unique solution
$Z_0 = \sum_{(i,j) \in \ZZ^2} a_{i,j} x^iy^j \in (p) \subset \ZZ_q\langle \ZZ^2 \rangle$ to
the equation $H(Z) = 0$. Moreover,
if $m \in \NN$ and $(r,s) \in \ZZ^2$ are such that
$(r,s) \notin m \Gamma$, then $\ord_p a_{r,s} \geq \frac{m}{2(a+b)}$.
\end{lemma}

\begin{remark} Note that we implicitly force $\Gamma$ to contain the origin: this follows from conditions 1 and 3.
If $\Gamma = \{ (0,0) \}$, Lemma~\ref{hensel} is just Hensel's lemma over $\ZZ_q$.
Finally, remark that if $(r,s)$ is not contained in any multiple of $\Gamma$, the above lemma implies
that $a_{r,s}$ equals $0$.
\end{remark}

\bproof
The existence and uniqueness of $Z_0$ follow immediately from Hensel's lemma, applied over $\ZZ_q \langle \ZZ^2 \rangle$.
Therefore we only need to
prove the convergence bound. Let $(r,s) \in \ZZ^2$ and $m \in \NN$ be such that $(r,s) \notin m\Gamma$.
Then there exists an edge spanning a line $eX + fY = c$ ($e,f,c \in \ZZ$), where $\Gamma \subset \{ (i,j) \in \ZZ^2 \ | \ ei + fj \leq c \}$,
such that $er + fs > mc$. Using a transformation of variables of
the type used in Lemma~\ref{fundamental}, we may assume that $e=0, f=1$ and $c \geq 0$. Thus $s>mc$.

Now, replace in $H(Z)$ all occurrences of $y^{-1}$ with a new variable $t$. We get
\[ H_\text{repl}(Z) = \sum h_{k, \text{repl}} (x,y,t) Z^k \in \ZZ_q[x^{\pm 1},y,t][Z] \]
with $\deg_y h_{k, \text{repl}} \leq (ak+b)c$.
Note that the conditions for Hensel's lemma are still satisfied. So there exists
a unique
\[ Z_{0, \text{repl}} = \sum_{(i,j,k) \in \ZZ \times \NN^2} b_{i,j,k} x^iy^jt^k \quad \in (p) \subset \ZZ_q \langle x^{\pm 1}, y, t \rangle\]
satisfying $H_\text{repl}(Z_{0, \text{repl}}) = 0$. If we substitute $y^{-1}$ for $t$, we get
precisely $Z_0$, due to the uniqueness statement in Hensel's lemma. Henceforth
\begin{equation} \label{frombtoa}
   a_{r,s} = \sum_{j - k = s} b_{r,j,k}.
\end{equation}
Let $K$ be a suitably ramified extension of $\QQ_q$ and denote by $R$ its valuation ring. Consider
\[ H'_\text{repl}(Z') = \sum p^{\mu_2(k-1)}h_k(x,p^{-\mu_1}y',t)Z'^k \quad \in K[x^{\pm 1},y',t][Z'] \]
obtained from $H_\text{repl}(Z)$ by substituting $y \leftarrow p^{-\mu_1}y'$, $Z \leftarrow p^{\mu_2}Z'$
and multiplying everything with $p^{-\mu_2}$. Here $\mu_1, \mu_2$ are positive rational numbers to be determined later.
We know that if
\begin{equation} \label{condit1}
  \mu_2 + j\mu_1 < 1 \quad \forall j \leq bc,
\end{equation}
\begin{equation} \label{condit2}
  j\mu_1 < 1 \quad \forall j \leq (a+b)c,
\end{equation}
and
\begin{equation} \label{condit3}
  (k-1)\mu_2 \geq j\mu_1 \quad \forall j \leq (ak+b)c
\end{equation}
for $k=2, \dots, \deg H$, then $H'_\text{repl}$ has integral coefficients and
$H'_\text{repl}(0) \equiv 0$ and $\frac{dH'_\text{repl}}{dZ'}(0) \equiv 1$ mod $P$.
Here $P$ is the maximal ideal of $R$. In that case, Hensel's lemma implies that there
is a unique $Z'_{0, \text{repl}} \in P \cdot R \langle x^{\pm 1},y',t \rangle$ such that
$H'_\text{repl}(Z'_{0, \text{repl}}) = 0$. Write
\[ Z_{0, \text{repl}}' = \sum_{(i,j,k) \in \ZZ \times \NN^2} b'_{i,j,k} x^iy'^jt^k\]
and perform reverse substitution to obtain that
\[ \sum_{(i,j,k) \in \ZZ \times \NN^2} p^{\mu_2}p^{j\mu_1}b'_{i,j,k}x^iy^jt^k \]
is a solution to $H_\text{repl}(Z) = 0$ in $P \cdot R \langle x^{\pm 1},y,t \rangle$. Again using
the uniqueness statement in Hensel's lemma we conclude that this is precisely $Z_{0, \text{repl}}$. As a consequence
\[ \ord_p b_{i,j,k} \geq j \mu_1 + \mu_2.\]
Using (\ref{frombtoa}) we find that $\ord_p a_{r,s} \geq s \mu_1 + \mu_2 > mc\mu_1 + \mu_2$. This gives the desired result, since
we can take $\mu_2 = \frac{2(a+b)c - bc}{2(a+b)c + \varepsilon}$ and $\mu_1 = \frac{1}{2(a+b)c + \varepsilon}$ for
any $\varepsilon \in \QQ_{>0}$.
\eproof

We are now ready to describe the construction of $\mathcal{F}_p$. In doing so, we will
systematically make a notational distinction between power series $g$ and the cosets $[g]$ (modulo $f$) they
represent (something which is usually not done in order to simplify notation).
Throughout, the assumptions about $\Gamma$ made at the beginning of Section~\ref{cohomsec} should be kept in
mind\footnote{In fact, for this section it suffices that $\Gamma$ contains the origin (not necessarily
as an interior point).}.

We will use a technique that was first described in~\cite{DVCab}. Suppose we can find
a $Z_0 \in \ZZ_q \langle \ZZ^2 \rangle^\dagger$ and polynomials $\delta_x, \delta_y \in \ZZ_q[\ZZ^2]$ such that
\[ [ f^\sigma(x^p(1 + \delta_x Z_0), y^p(1 + \delta_y Z_0))] = [0] \quad \text{in $A^\dagger$}, \]
where $f^\sigma$ is obtained from $f$ by applying Frobenius substitution\footnote{By Frobenius
substitution we mean the map $\ZZ_q \rightarrow \ZZ_q : \sum_{i=0}^\infty \pi_i p^i
\mapsto \sum_{i=0}^\infty \pi_i^p p^i$, where the $\pi_i$ are Teichm\"uller representatives.}
to the coefficients.
Then
\[ \mathcal{F}_p : A^\dagger \rightarrow A^\dagger :
\left\{ \begin{array}{lll}
\text{$[x]$} & \mapsto & \text{$[x^p(1 + \delta_xZ_0)]$} \\
\text{$[y]$} & \mapsto & \text{$[y^p(1 + \delta_yZ_0)]$} \\
\end{array} \right. \]
(acting on $\ZZ_q$ by Frobenius substitution and extended by linearity and continuity)
is a well-defined $\ZZ_p$-algebra morphism that lifts $\overline{\mathcal{F}}_p$.

Take $\overline{\beta}, \overline{\beta}_x, \overline{\beta}_y \in \FF_q[\ZZ^2]$ with support in $2 \Gamma$
for which
\[ 1 = \overline{\beta} \, \overline{f} + \overline{\beta}_x x \frac{\partial \overline{f}}{\partial x}
+ \overline{\beta}_y y \frac{\partial \overline{f}}{\partial y} \]
(this is possible due to Corollary~\ref{nullbound}). Let $\delta, \delta_x, \delta_y$ be arbitrary
Newton polytope preserving lifts of $\overline{\beta}^p, \overline{\beta}_x^p$ resp. $\overline{\beta}_y^p$.
Then clearly $\Gamma(\delta), \Gamma(\delta_x), \Gamma(\delta_y) \subset 2p\Gamma$.

Now let $x^ay^b$ be any monomial such that $g(x,y) = x^ay^b f(x,y)$ has support in $\NN^2$
and define $G(Z) = x^{-pa}y^{-pb}g^\sigma (x^p(1 + \delta_xZ), y^p(1 + \delta_yZ)) \in \ZZ_q[\ZZ^2][Z]$, where
$g^\sigma$ is again obtained from $g$ by applying Frobenius substitution to the coefficients. Since
\[ G(0) \equiv f^p \quad \text{and} \quad \frac{dG}{dZ}(0) \equiv 1 + (a \delta_x + b \delta_y - \delta)f^p \quad \text{mod $p$} \]
we see that $[G(Z)]=[0]$ has a unique solution $[Z_1]$ that is congruent to $0$ mod $p$
in the Henselian ring $A^\dagger$. However, Hensel's lemma does not provide any information on the
convergence rate of $Z_1$ (or any other representant of $[Z_1]$). To solve this problem, define
\[ H(Z) = G(Z) - (a \delta_x + b \delta_y - \delta)f^pZ - f^p.\]
Then clearly $[G(Z)] = [H(Z)]$, but now the conditions of Hensel's lemma are satisfied over
the base ring, so that there exists a unique $Z_0 \in (p) \subset \ZZ_q \langle \ZZ^2 \rangle$
for which $H(Z_0) = 0$. We have that $[Z_0] = [Z_1]$. Note that if we expand
\[ H(Z) = \sum_{k=0}^{\deg H} h_k(x,y) Z^k,\]
one easily checks that
\begin{equation} \label{coefficientsofH}
\Gamma(h_k) \subset (2k+1) p \Gamma.
\end{equation}
Therefore, we can apply Lemma~\ref{hensel}
and conclude that $Z_0 = \sum_{(i,j) \in \ZZ^2} a_{i,j} x^iy^j$ where the $a_{i,j}$ satisfy:
\begin{equation} \label{boundforrootsofH}
   \forall i,j \in \ZZ, m \in \NN : \ (i,j) \notin m \Gamma \, \Rightarrow \ord_p a_{i,j} \geq \frac{m}{6p}.
\end{equation}

Our next step is to investigate what the convergence rate of $Z_0$ tells us about the convergence rate
of $Z_x = 1 + \delta_xZ_0$ and $Z_y = 1 + \delta_yZ_0$.
Write $Z_x = \sum_{(i,j) \in \ZZ^2} b_{i,j} x^iy^j$. We claim that
\[ \forall i,j \in \ZZ, m \in \NN : \ (i,j) \notin m \Gamma \, \Rightarrow \ord_p b_{i,j} \geq \frac{m}{8p}.\]
Indeed, since $Z_0 \equiv 0$ mod $p$, this statement is definitely true for $m < 8p$.
If $m \geq 8p$, then $\frac{m-2p}{6p} \geq \frac{m}{8p}$. Now suppose $(i,j) \notin m\Gamma$. Write
$\delta_x = \sum_{(i,j) \in 2p\Gamma} d_{i,j}x^iy^j$. We know that
\[ b_{i,j} = \sum_{k + r = i, \, \ell + s = j} d_{k, \ell} a_{r,s} \]
and since $(k, \ell) \in 2p\Gamma$, we know that all $(r,s)$ appearing in the above expansion are not contained
in $(m-2p)\Gamma$. Therefore
\[ \ord_p b_{i,j} \geq \frac{m-2p}{6p} \geq \frac{m}{8p}.\]
These observations allow us to state the main result of this section.

\begin{theorem} \label{thmlift}
There exist units $Z_x, Z_y \in \ZZ_q \langle \ZZ^2 \rangle^\dagger$ such that
\[ \mathcal{F}_p : A^\dagger \rightarrow A^\dagger :
\left\{ \begin{array}{lll}
\text{$[x]$} & \mapsto & \text{$[x^pZ_x]$} \\
\text{$[y]$} & \mapsto & \text{$[y^pZ_y]$} \\
\end{array} \right.\]
(extended by linearity and continuity and acting on $\ZZ_q$ by Frobenius substitution) is a well-defined $\ZZ_p$-algebra
morphism that lifts $\overline{\mathcal{F}}_p$. Moreover $Z_x, Z_y, Z_x^{-1}, Z_y^{-1}$ satisfy
the following convergence criterion: if $(i,j) \in \ZZ^2, m \in \NN$ are such that $(i,j) \notin m \Gamma$, then the
coefficient of $x^iy^j$ has $p$-order $> \frac{m}{9p}$.
\end{theorem}

\bproof
It only remains to show that $Z_x^{-1}$ and $Z_y^{-1}$ satisfy the convergence criterion. This can be done
as in the proof of Lemma~\ref{hensel}. We refer to \cite{Wouterthesis} for a detailed proof.
\eproof

\begin{remark}
The larger denominator ($9p$ instead of $8p$) is a small price we
have to pay during inversion, but it also allows us to write down a strict inequality
($>$ instead of $\geq$). In this form, the convergence criterion
is closed under multiplication, i.e.
\begin{equation} \label{ring}
  \left\{ \left. \sum_{(i,j) \in \ZZ^2} a_{i,j}x^iy^j \in \ZZ_q\langle\ZZ^2\rangle \, \right| \,
\forall m \in \NN, (i,j) \in \ZZ^2 : (i,j) \notin m\Gamma \Rightarrow \ord_pa_{i,j} > \frac{m}{9p} \right\}
\end{equation}
is a ring. We will use this in Section~\ref{complsec}.
\end{remark}


\section{Reduction Algorithm}\label{redsec}

Throughout this section, $\Gamma$ should again satisfy the assumptions made at the beginning of Section~\ref{cohomsec}.
An important step in our algorithm (see Section~\ref{complsec}) is to reduce Laurent polynomials modulo the operator
$D$ defined in Section~\ref{cohomsec}. Below we describe a procedure
that solves this problem and prove that,
after multiplying with a small power of $p$, the reduction process is
entirely integral, which enables us to tightly bound the loss of precision.

First, we need a few more theoretical results. Let $\{t_1, \dots, t_r\}$ be the edges of $\Gamma$,
let $\{T_1, \dots, T_r\} \subset X_\Gamma$ be
the corresponding $\QQ_q$-tori and let $\{\overline{T}_1, \dots, \overline{T}_r\}
\subset \overline{X}_\Gamma$ be the corresponding $\FF_q$-tori. The reductions mod $p$ of the points
in $T_k \cap C$ are precisely the points of $\overline{T}_k \cap \overline{C}$.
For every $k=1, \dots, r$, we can
find $c_k,b_k \in \ZZ$ such that
\begin{equation} \label{locpartype}
x^{c_k}y^{b_k}
\end{equation}
defines a local parameter, at both $P$ and $\overline{P}$, for each $P \in T_k \cap C$.
Indeed, these assertions follow from the proof of Lemma~\ref{fundamental}: $c_k, b_k$
depend only on the
geometry of $\Gamma$. If in what follows we say `\emph{local parameter over $\ZZ_q$}', actually
any $t \in \ZZ_q [\ZZ^2]/(f)$ for which both $t$ and its reduction mod $p$ are local parameters
at $P$ resp. $\overline{P}$ will work.

Below, let
$\QQ_q^\text{ur} \subset \overline{\QQ}_q$ denote the maximal unramified extension of $\QQ_q$ and
let $\ZZ_q^\text{ur}$ be its valuation ring. Note that
all places $P \in C \setminus \TT^2_{\QQ_q}$ are defined over $\QQ_q^\text{ur}$.

\begin{definition} \mbox{} \\[-0.5cm]
\begin{enumerate}
\item Let $L^{(0)} = \ZZ_q^\emph{ur}[\ZZ^2]$, then for any set $S$
of Laurent polynomials, define $S^{(0)} = S \cap L^{(0)}$.
\item Let $L^{(1)}$ be the subset of $L^{(0)}$ consisting of those $h$ for which the following
holds. For every $P \in C \setminus \TT^2_{\Qq}$, take a local parameter $t$ over $\ZZ_q$. Then the
condition is that
\[ \frac{t}{dt} \Lambda(h) = \sum_{i = v}^\infty a_i t^i \quad (a_i \in \ZZ^\emph{ur}_q)\]
satisfies $\ord_p a_i \geq \ord_p i$ (alternative notation: $i | a_i$) for all $i<0$.
For any set $S$ of Laurent polynomials let
$S^{(1)} = S \cap L^{(1)}$.
\end{enumerate}
\end{definition}

Again we remark that the above definitions are vulnerable to notational abuses. For instance,
if $S$ consists of \emph{cosets} of Laurent polynomials, then $S^{(0)}$ consists of those Laurent
polynomials having a representant in $L^{(0)}$, and so on.

The set $L^{(1)}$ appears naturally\footnote{In \cite[Proposition 5.3.1]{Edixhoven},
Edixhoven uses a similar set.} when
we apply the operator $D = xy\left(\frac{\partial f}{\partial y} \frac{\partial}{\partial x}
- \frac{\partial f}{\partial x} \frac{\partial}{\partial y} \right)$ that was introduced in Section~\ref{cohomsec}
to an element in $L^{(0)}$.

\begin{lemma}
If $h \in L^{(0)}$, then $D h  \in L^{(1)}$.
\end{lemma}
\bproof
Let $P \in C \setminus \TT^2_{\Qq}$ and let $t$ be a local parameter over $\ZZ_q$ at $P$.  By the definition of $D$, we have
\[ \frac{t}{dt} \Lambda(D h) = \frac{t}{dt} dh  \, . \]
Write $h = \sum_{i = v}^{\infty} b_i t^i$, then clearly $\frac{t}{dt} dh = \sum_{i = v}^{\infty} i b_i t^i$,
which proves the claim.
\eproof

\begin{lemma} \label{strongRRSpace}
Let $D$ be a divisor on $C$ which is defined over $\QQ_q$ and which has support in $C \setminus \TT_{\QQ_q}^2$. Then $\mL^{(0)}(D)$ is
free and finitely generated over $\ZZ_q$.
\end{lemma}
\bproof
We first prove that the following `strong' version of Theorem~\ref{RRSpace} holds: \emph{for every
$m \in \mathbb{N}_0$, the module $\mL^{(0)}(mD_C)$ is precisely given by $L^{(0)}_{m\Gamma}$.}
Take an element of $\mL^{(0)}_{m\Gamma}$, represented by some
$h \in \ZZ_q[\ZZ^2]$. By Theorem~\ref{RRSpace}, there is an $\alpha \in \QQ_q[\ZZ^2]$
such that $h + \alpha f$ has support in $m\Gamma$. Write $\alpha = \alpha_1 + \alpha_2$,
where all coefficients of $\alpha_1$ are integral and all coefficients of $\alpha_2$ are
non-integral. We claim that $h + \alpha_1 f$ has support in $m\Gamma$. Indeed, suppose
this were not true, then $\alpha_2 f$ has a non-zero term with support outside
$m\Gamma$. This implies that $\alpha_2$ has a non-zero term with support outside
$(m-1)\Gamma$. Let $a_{ij}x^iy^j$ be such a term. Then $\Gamma$ has an edge
spanning a line $dX + eY = c$ (with $\Gamma \subset \{ (r,s) \, | \, dr + es \leq c \}$)
such that $di + ej > (m-1)c$. Consider the following monomial order:
\[ \begin{array}{lrl} x^ry^s \prec x^ky^\ell & \text{if} & dr + es < dk + e\ell \\
                                             & \text{or if} & dr + es = dk + e\ell \ \text{and} \ r < k \\
                                             & \text{or if} & dr + es = dk + e\ell, \ r = k \ \text{and} \ s < \ell \\
\end{array} \]
(where the last line is only of use if $e=0$).
We may suppose that $x^iy^j$ is maximal with respect to $\prec$.
Take the term $b_{rs}x^ry^s$ of $f$ that is maximal with respect to $\prec$ (in particular, $dr + es = c$).
Then $a_{ij}b_{rs}x^{i+r}y^{j+s}$ is a term of $\alpha_2f$ with support outside $m\Gamma$. Because
$h + \alpha_1f + \alpha_2f$ has support in $m\Gamma$ and $h + \alpha_1f \in \ZZ_q[\ZZ^2]$, this implies
that $a_{ij}b_{rs}$ is integral. But this is impossible, since $a_{ij}$ is non-integral and $b_{rs}$ is
a $p$-adic unit.

Now since $\mL^{(0)}(D) \subset \mL^{(0)}(mD_C)$ for some big enough $m \in \NN_0$, and since the latter
is finitely generated, we have that $\mL^{(0)}(D)$
is finitely generated as well. This follows from a well-known theorem on modules over Noetherian rings.
But it is also well-known that every finitely generated and torsion-free module over
a principal ideal domain is free, which concludes the proof.
\eproof

For the following two lemmata, fix a point $P \in C \setminus \TT^2_{\QQ_q}$ and let
$\QQ_{q^s} \supset \QQ_q$ be its field of definition. Denote the valuation ring with
$\ZZ_{q^s}$ and the residue field with $\FF_{q^s}$. Write $\text{Gal}(\QQ_{q^s}, \QQ_q) = \{ \sigma_1, \sigma_2, \dots, \sigma_s \}$
with $\sigma_1 = \text{id}_{\QQ_{q^s}}$. Let $\textbf{P}$ be the divisor $\sum_{i=1}^s P^{\sigma_i}$.
Note that if $t$ is a local parameter at $P$ over $\ZZ_q$, then it is a local parameter
at any $P^{\sigma_i}$ over $\ZZ_q$.

\begin{lemma}\label{existh}
Let $E$ be an effective divisor on $C$ which is defined over $\QQ_q$
and whose support is contained in $C \setminus \TT^2_{\Qq}$.
Assume that $\deg E > 2g - 2$. Then there exists an
$h \in \mL^{(0)}_{\QQ_q^s}(E + P)$ such that:
\begin{enumerate}
\item $h$ has a pole at $P$ of multiplicity $\ord_P(E) + 1$.
\item Let $t$ be a local parameter over $\ZZ_q$ at $P$. Then $h$ has an expansion
 $\sum_{i = v}^\infty a_i t^i$,
with all $a_i \in \ZZ_{q^s}$ and $a_v$ a unit in $\ZZ_{q^s}$.
\end{enumerate}
\end{lemma}
\bproof
Consider the following diagram where the vertical arrows are the natural reduction
modulo $p$ maps:
\[
\xymatrix{
{\mL}^{(0)}_{\QQ_{q^s}}(E) \, \ar@{->}^{\subset \quad }[r] \ar@{->>}[d] & {\mL}^{(0)}_{\QQ_{q^s}}(E+P) \ar@{->>}[d] \\
{\mL}_{\oC, \FF_{q^s}}(\overline{E}) \ar@{->}^{\subsetneq \quad }[r] & {\mL}_{\oC, \FF_{q^s}}(\overline{E} + \overline{P}).
}
\]
The vertical maps are surjective, since after tensoring with $\Fq$ they become
clearly injective and hence surjective since both have the same dimension
by Riemann-Roch (here we used the foregoing lemma).
Let $\overline{h} \in {\mL}_{\oC, \FF_{q^s}}(\overline{E} + \overline{P}) \setminus {\mL}_{\oC, \FF_{q^s}}(\overline{E})$
and choose $h \in {\mL}^{(0)}_{\QQ_{q^s}}(E+P)$ such that $h$ reduces to $\overline{h}$ mod $p$.
\eproof

An important feature of the foregoing lemma is the following: if we replace
$P$ by $P^{\sigma_i}$ for some $\sigma_i \in \text{Gal}(\QQ_{q^s}, \QQ_q)$,
then $h^{\sigma_i} \in \mathcal{L}^{(0)}_{\QQ_{q^s}}(E + P^{\sigma_i})$
again satisfies the above conditions. Indeed, $\sigma_i(a_v)$ is a unit
in $\ZZ_{q^s}$.

\begin{lemma} \label{toberepeated}
Let $E$ be an effective divisor on $C$ which is defined over $\QQ_q$ and whose
support is contained in $C \setminus \TT^2_{\Qq}$.  Suppose that $\deg E > 2g - 2$, then the map
\[ \mL^{(0)} (E + \emph{\textbf{P}}) \stackrel{D}{\longrightarrow} \frac{\mL^{(1)}(E + D_C + \emph{\textbf{P}})}{\mL^{(1)}(E + D_C)} \]
is surjective.
\end{lemma}
\bproof
Let $h \in \mL^{(1)}(E + D_C + \textbf{P})  \setminus \mL^{(1)}(E + D_C)$. By Corollary~\ref{canonicaldivisor} we have
$\Div \Lambda(h) = \Div h + D_C - W_C$.  Let $t$ be a local parameter over $\ZZ_q$ at $P$, then
\[
\begin{split} \ord_P \left( \frac{t \Lambda(h)}{dt} \right) & = \ord_P ( h ) + \ord_P(D_C)  \\
& = - \ord_P(E + D_C + \textbf{P}) + \ord_P(D_C) \\
& = - \ord_P(E) - 1 = -n \,,
\end{split}
\]
with $n = \ord_P(E + \textbf{P})$.  Therefore we have a local expansion
\[ \frac{t \Lambda(h)}{dt} = b_0 t^{-n} + b_1 t^{-n+1} + \cdots  \, , \]
at $P$, with $n | b_0$. Note that the expansions at
the conjugate places $P^{\sigma_i}$ are given by
\[ \frac{t \Lambda(h)}{dt} = \sigma_i(b_0) t^{-n} + \sigma_i(b_1) t^{-n+1} + \cdots . \]
Using Lemma~\ref{existh}
we find an $h_0 \in \mL^{(0)}_{\QQ_{q^s}}(E + P)$ with
power series expansion at $P$:
\[ h_0 = a_0 t^{-n} + a_1 t^{-n+1} + \cdots  \, , \]
and with $a_0$ a $p$-adic unit.
Define
\[ h_1 = h + \sum_{i=1}^s \frac{\sigma_i(b_0)}{n \sigma_i(a_0)} D(h_0^{\sigma_i}) = h + D\left(\text{Tr}\left( \frac{b_0}{n a_0} h_0 \right) \right) \in \mL^{(1)}(E + D_C + \textbf{P}),\]
then we have the following expansion at $P$:
\[ \frac{t \Lambda(h_1)}{dt}  = \frac{t \Lambda(h)}{dt} + \frac{b_0}{n a_0} \frac{t d h_0}{dt} = 0 \cdot t^{-n} + \cdots \, ,\]
and thus $\ord_P \left( \frac{t \Lambda(h_1)}{dt} \right)  \geq -n + 1$. Similarly,
the pole orders at all conjugate places $P^{\sigma_i}$ are reduced by at least $1$.
Note that
\[ \ord_{P^{\sigma_i}} \left( \frac{t \Lambda(h_1)}{dt} \right) = \ord_{P^{\sigma_i}} (h_1) + \ord_{P^{\sigma_i}}(D_C) \, , \]
since $\Div \Lambda(h_1) = \Div h_1 + D_C - W_C$.  Hence we see that
$\ord_{P^{\sigma_i}}(h_1) \geq -n + 1 - \ord_{P^{\sigma_i}}(D_C) = 1 - \ord_{P^{\sigma_i}}(E + D_C + \textbf{P})$, thus $h_1 \in \mL^{(1)}(E + D_C)$
which finishes the proof.
\eproof

A repeated application of the above lemma gives the following result.
\begin{corollary} \label{reductioncoro}
Let $E$ be an effective divisor which is defined over $\QQ_q$ and
whose support is contained in $C \setminus \TT^2_{\Qq}$, then the map
\[ \mL^{(0)} (D_C + E) \stackrel{D}{\longrightarrow} \frac{\mL^{(1)}(2 D_C + E)}{\mL^{(1)}(2 D_C)} \]
is surjective.
\end{corollary}

The above corollary can be turned into a reduction algorithm and also provides a sharp bound
for the loss of precision incurred during reduction.  Indeed, since the Newton polytope
$\Gamma$ contains the origin as an interior point, any Laurent polynomial $h \in \Zq[\ZZ^2]$ will be contained in an
$L^{(0)}_{m \Gamma}$ with $m \in \NN_0$ big enough.  Let
\[ \varepsilon =  \ceil{ \log_p \max \{ - \ord_P(h) \}_{P \in C \setminus \TT^2_{\Qq}}}  \, , \]
then clearly $p^\varepsilon h \in L^{(1)}_{m \Gamma}$. So we can as well assume that $h \in L^{(1)}_{m\Gamma}$.
By Theorem~\ref{RRSpace}, we have $\mL(m D_C) = L_{m \Gamma}$ and applying Corollary~\ref{reductioncoro}
with $E = (m-2) D_C$ (we can assume that $m > 2$, since otherwise no reduction is necessary), shows
that there exists a $g \in L^{(0)}_{(m-1)\Gamma}$ such that $h_r = h - D(g) \in {\mL^{(1)}(2 D_C)}$.
Note that after multiplication with $p^\varepsilon$ the entire reduction process is integral,
so if we want to recover the result $h_r$ modulo $p^N$, we need to compute $h$ modulo
$p^{N + \varepsilon}$.  To finalize the computation, we need to express $h_r$
on a basis for $H^{1}_{DR}(C \cap T^2_{\Qq})$, which could cause a further loss of precision,
depending on the basis chosen.  But clearly, as long a we choose a `$\Zq$-module basis' for
$H^{1}_{DR}(C \cap T^2_{\Qq})$, no further loss of precision will occur. More precisely, we mean
the following. Consider the module
\[ M_H = \frac{\mL^{(0)}(2 D_C)}{D (\mL (D_C)) \cap L^{(0)} }  \, , \]
then $M_H$ is a free $\Zq$-module since it is finitely generated and torsion-free.
Therefore, any $\Zq$-basis for $M_H$ forms a suitable basis for $H^{1}_{DR}(C \cap T^2_{\Qq})$,
such that in the final reduction step, no further loss of precision is incurred.

In the above description, we used any representant for an element of the coordinate
ring of $C$; in practice however, we would like to work with a unique representant.
Given the Newton polytope $\Gamma$ of $f$, there are many possibilities to choose
a suitable basis $B$ for $\Qq[\ZZ^2]/(f)$.  The assumptions about $\Gamma$ made in
Section~\ref{cohomsec} already led to the following natural choice
\[ B = \{ x^k y^l \ | \ k, l \in \ZZ, d_b \leq l < d_t \} \, , \]
with $(c_t, d_t)$ (resp.\ $(c_b, d_b)$ ) the unique highest (resp.\ lowest) point
of $\Gamma$.

Let $S_{[m_1,m_2]}$ with $m_1 < m_2$ denote the set of Laurent polynomials with support
in the rectangle $[m_1, m_2] \times [d_b, d_t - 1]$, then the reduction process proceeds
in two phases: the first phase reduces terms in $S_{[0,m]}$ with $m \in \NN_0$ and
the second phase reduces terms in $S_{[-m,0]}$ with $m \in \NN_0$.  Since both
phases are so similar, we will focus mainly on the first phase and briefly
mention the changes for the second phase.


\subsubsection{Phase 1:}

Any element $h \in S^{(0)}_{[0,m]}$ can be forced into
$S^{(1)}_{[0,m]}$ by multiplying it with $p^\varepsilon$ where
\[ \varepsilon = \ceil{ \log_p (m M_x + \Delta) } \]
with $M_x = \max \{ - \ord_P (x) \}_{P \in C \setminus
\TT^2_{\Qq}}$ and $\Delta = \max \{ - \ord_P (y^{d_t-1}) , -
\ord_P (y^{d_b}) \}_{P \in C \setminus \TT^2_{\Qq}}$. If we now
want to apply Corollary~\ref{reductioncoro} to an element $h \in
S^{(1)}_{[0,m]}$, we need to find a divisor $E$ over $\QQ_q$ such
that $S^{(1)}_{[0,m]} \subset \mL^{(1)}(2 D_C + E)$. Then by
Corollary~\ref{reductioncoro} there exists a $g \in \mL^{(0)} (D_C
+ E)$ such that $h - D(g) \in \mL^{(1)}(2D_C)$.  In practice
however, we do not want to work with explicit Riemann-Roch spaces;
as such we want to find a divisor $E$ (depending on $m$) and
constants $c_1, c_2 \in \ZZ$ (independent of $m$) such that
\[ S_{[0,m]} \subset \mL(2D_C + E)  \quad \text{ and } \quad \mL(D_C + E) \subset S_{[c_1, m + c_2]} \, . \]
The reduction algorithm then becomes very simple indeed: to reduce
$h \in S^{(1)}_{[0,m]}$, we only need to find a $g \in
S^{(0)}_{[c_1, m + c_2]}$ such that $h - D(g) \in
\mL^{(1)}(2D_C)$, using linear algebra.

Recall that the divisor of any function $h \in \Qq(C)$ can be written as
the difference of the zero divisor and the pole divisor, i.e.\
$\Div(h) = \Div_0(h) - \Div_\infty(h)$, $\Div_0(h) \geq 0$, $\Div_\infty(h) \geq 0$
and $\Supp (\Div_0(h)) \cap \Supp (\Div_\infty (h)) = \emptyset$.  Furthermore,
two trivial observations are that $h \in \mL(\Div_\infty(h))$ and
$\Div_\infty(h^{-1}) = \Div_0(h)$.  Consider the divisor
\[ E_m =  -d_b \Div_0 (y) + (d_t - 1) \Div_\infty (y)  + m \Div_\infty (x) \, \]
then $E_m \geq 0$ and $S_{[0,m]} \subset \mL(E_m) \subset \mL(2 D_C + E_m)$,
so we can apply Corollary~\ref{reductioncoro} with $E = E_m$.
Note that $E_m$ is indeed defined over $\QQ_q$.

\begin{remark} \label{remEm}
It is clear that the choice for $E_m$ is not entirely optimal, since we could subtract
the contributions in $2 D_C$ and still obtain the above inclusion.  The most important
simplification in practice is that $2 \Gamma$ is `likely' to contain the
interval $[d_b, d_t-1]$ on the $y$-axis and then $E_m$ can be simply taken to be $m \Div_\infty(x)$.
However, in general this need not be the case.
\end{remark}

To determine the constants $c_1$ and $c_2$ we first prove the following lemma.
\begin{lemma}\label{rrsetinclusion}
Let $E$ be a divisor on $C$ which is defined over $\QQ_q$ and with $\deg E > 2g - 2$, and let $h \in \Qq(C)$ be
a function on $C$. Then for any $m \in \NN_0$ the following map is an isomorphism:
\[ \frac{\mL(E + \Div_\infty(h))}{\mL(E)} \ \xrightarrow{\cdot h^{m-1}} \
\frac{\mL(E + m \Div_\infty(h))}{\mL(E + (m-1) \Div_\infty(h))} \, .
\]
\end{lemma}
\bproof
Since $\deg E > 2g - 2$ and $\Div_\infty (h) \geq 0$, the Riemann-Roch theorem
implies that the dimensions of both vector spaces are equal to $\deg \Div_\infty(h)$,
so it suffices to prove injectivity.  Let $g \in \mL(E + \Div_\infty(h))$
and assume that $h^{m-1} g \in \mL(E + (m-1) \Div_\infty(h))$, i.e.\
\[ (m-1) \Div (h) + \Div(g) \geq - E - (m-1) \Div_\infty(h) \, , \]
which implies that $\Div(g) \geq - E - (m-1) \Div_0 (h)$.  Since $g \in \mL(E + \Div_\infty(h))$,
i.e.\ $\Div(g) \geq - E - \Div_\infty (h)$ and the supports of $\Div_0(h)$ and $\Div_\infty(h)$
are disjoint, we conclude $\Div(g) \geq - E$ or $g \in \mL(E)$.
\eproof

In what follows, we will use the abbreviation $E_y = -d_b \Div_0 (y) + (d_t - 1) \Div_\infty (y)$,
so $E_m = E_y + m\Div_\infty(x)$.
Choose integers $\kappa_1 \leq 0$ and $\kappa_2 \geq 0$ such that
$\mL^{(0)}(D_C + E_y + \Div_\infty(x) + \Div_0(x)) \subset S_{[\kappa_1, \kappa_2]}$. In particular,
$\mL^{(0)}(D_C + E_1) \subset S_{[\kappa_1, \kappa_2]}$. This can then be generalized to the following.

\begin{corollary} \label{defkappa}
$\mL(D_C + E_m) \subset S_{[\kappa_1, m - 1 + \kappa_2 ]}$.
\end{corollary}
\bproof
Apply Lemma~\ref{rrsetinclusion} with $E = D_C + E_y$ and $h = x$.
\eproof

Thus, given $h \in S^{(1)}_{[0,m]}$ we find $g \in S^{(0)}_{[\kappa_1, m - 1 + \kappa_2 ]}$ such that
$h - D(g) \in \mL^{(1)}(2 D_C)$ using linear algebra over $\Zq$.  However, for big $m$ the
linear systems involved get quite large, so we compute $g$ in several steps: let $h_0 = h$
and choose a constant $c \in \NN_0$, then in step $1 \leq i \leq t$ (where $t$ will be determined
later) we compute a $g_i$ such that
\[ h_i = h_{i-1} - D(g_i) \in S^{(1)}_{[0, m - ic]} \, . \]
In the last step, i.e.\ step $t+1$ we
find a $g_{t+1} \in S^{(0)}_{[\kappa_1, m - tc - 1 + \kappa_2 ]}$ such
that
\[ h_{t+1} = h_{t} - D(g_{t+1}) \in \mL^{(1)}(2 D_C) \, .  \]
We postpone this last step until after Phase 2, since
it is better to treat the last steps of both phases at once.
To determine which monomials appear in the $g_i$ for $1 \leq i \leq t$
we prove the following lemma.

\begin{lemma}\label{redboundlem}
If $m \in \NN_0, k \in \ZZ$ with $d_b \leq k < d_t$, then $D(x^m
y^k) \in S^{(1)}_{[\kappa_1 + m - 1, \kappa_2 + m - 1]}$.
\end{lemma}
\bproof
By definition of $D$ we have $D(x^m y^k) = x^m y^k (m y f_y - k x f_x)$.  Note
that the support of $g = m y f_y - k x f_x$ is contained in $\Gamma$ and thus
$g \in \mL(D_C)$.  Furthermore, by definition of $E_y$ we have $y^k \in \mL(E_y)$.
Therefore, by definition of $\kappa_1$ and $\kappa_2$ we conclude that
$D(x^m y^k) \in S^{(1)}_{[\kappa_1 + m - 1, \kappa_2 + m - 1]}$.
\eproof

The above lemma finalizes the description of the algorithm: in step $i$
it suffices to take $g_i$ in $S_{[a_i, b_i]}$ with
\[ a_i = m - ic - \kappa_2 + 2 \quad \text{ and } \quad b_i = m - (i-1)c + \kappa_2 - 1 \, , \]
and to work modulo $x^{m - ic}$.
There are two natural conditions that $t$ and $c$ should satisfy.
The first one is related to the fact that we want to work in $S_{[0,+\infty]}$ only. Therefore,
\[ a_t \geq - \kappa_1 + 1 \quad \text{ which is equivalent with } \quad t c \leq m + \kappa_1 - \kappa_2 + 1 \, . \]
The second condition keeps track of the fact that something which is
already in $\mL^{(1)}(2D_C)$ cannot be reduced anymore. Therefore, choose integers\footnote{The parameters $\kappa_1, \kappa_2$
and $\chi_1, \chi_2$ will be discussed more extensively in Section~\ref{complsec}.}
$\chi_1 \leq 0, \chi_2 \geq 0$
such that $\mL^{(1)}(2D_C) \subset L_{[2\chi_1, 2\chi_2]}$. It then suffices to impose
\[ tc \leq m - 2\chi_2.\]
The number of unknowns in the linear system of equations in step $i$ is precisely
the number of monomials in $S_{[a_i, b_i]}$, which equals $(d_t - d_b)(c + 2\kappa_2 - 2)$.
Note that this also appears as a natural upper bound for the number of terms in $D(S_{[a_i,b_i]})$ modulo $x^{m -ic}$,
so we obtain a system with as least as many unknowns as equations.



\subsubsection{Phase 2:}  Since the second phase is very similar to the first, we will only
briefly mention the main differences.  To force an element $h \in S^{(0)}_{[-m,0]}$
with $m \in \NN_0$ into $S^{(1)}_{[-m,0]}$, we need to multiply with
$p^\varepsilon$ where
\[ \varepsilon = \ceil{ \log_p (m M_{1/x} + \Delta) } \]
with $M_{1/x} = \max \{ - \ord_P (x^{-1}) \}_{P \in C \setminus \TT^2_{\Qq}}$ and $\Delta$
as before, so from now on assume that $h \in S^{(1)}_{[-m,0]}$.
The divisor $E_m$ now becomes $E_m = E_y + m \Div_\infty(x^{-1})$ and
applying Lemma~\ref{rrsetinclusion} with $h = x^{-1}$ shows
\[ \mL(D_C + E_y + m \Div_\infty (x^{-1})) \subset S_{[-m + 1 + \kappa_1, \kappa_2 ]} \, , \]
where $\kappa_1, \kappa_2$ are chosen as in Phase 1. In step $i$
we now compute a $g_i$ such that $h_i = h_{i-1} - D(g_i) \in S^{(1)}_{[-m + ic, 0]}$ for
some constant $c \in \NN_0$.  An analogue of Lemma~\ref{redboundlem}
(replace $S^{(1)}_{[\kappa_1 + m - 1, \kappa_2 + m - 1]}$ with $S^{(1)}_{[\kappa_1 - m + 1, \kappa_2 - m + 1]}$)
finally leads to
$g_i \in S^{(0)}_{[a_i, b_i]}$ with
\[  a_i = -m + (i-1)c + \kappa_1 + 1 \quad \text{ and } \quad b_i = -m + ic - \kappa_1 - 2 \, . \]
The number of steps $t$ is determined by the following inequalities:
\[ tc \leq m + \kappa_1 - \kappa_2 + 1 \quad \text{and} \quad tc \leq m + 2\chi_1.\]
The systems to be solved have $(d_t - d_b)(c - 2\kappa_1 - 2)$ unknowns, that are related by at most
the same number of equations.


\subsubsection{Step $t+1$:} During Phase 1 and Phase 2, we reduced a given polynomial
$h \in L^{(1)}$ modulo $D$ to obtain a polynomial $h_t \in S^{(1)}_{[-n_1,n_2]}$, where
$n_1 \in \mathbb{N}_0$ is roughly of size $\max\{ - 2\chi_1, \kappa_2 - \kappa_1 \}$ and
$n_2 \in \mathbb{N}_0$ is roughly of size $\max\{ 2\chi_2, \kappa_2 - \kappa_1 \}$. In this last
step, we reduce to a polynomial $h_{t+1} \in \mathcal{L}^{(1)}(2D_C)$ by brute force.
From Corollary~\ref{defkappa} (and its Phase 2 analogue) we
know that there is a $g_{t+1} \in S^{(0)}_{[-n_1 + 1 + \kappa_1, n_2 - 1 + \kappa_2]}$ such that
\[ h_t - D(g_{t+1}) \in \mathcal{L}^{(1)}(2D_C),\]
so we can compute $h_{t+1}$ by solving a system of at most $(d_t - d_b)(2(\kappa_2 - \kappa_1) + n_1 + n_2 - 3)$
equations in
\[ (d_t - d_b)(\kappa_2 - \kappa_1 + n_1 + n_2 - 1) + \#(2\Gamma \cap \ZZ^2)\]
unknowns. Here, the latter term equals $4\Vol(\Gamma) + \#(\partial \Gamma \cap \ZZ^2) + 1$
by Ehrhart's theorem \cite{Ehrhart}.


\subsubsection{Solving linear systems over $\ZZ_q$.}
Let $r,s \in \NN_0$ and consider a matrix $A \in \ZZ_q^{r \times s}$ and a vector $b \in \ZZ_q^r$. Let
$N \in \NN_0$ denote the $p$-adic precision up to which is to be computed. The aim is to
find an $\textbf{x} \in \ZZ_q^s$ such that $A \cdot \textbf{x} \equiv b$ mod $p^N$. Note that
this is slightly weaker than finding the reduction mod $p^N$ of an $\textbf{x} \in \ZZ_q^s$
such that $A \cdot \textbf{x} = b$ (exact equality over $\ZZ_q$), but only slightly: from Lemma~\ref{invariant}
below it follows that it suffices to increase the precision in order to solve this.

Using Gaussian elimination, where in each step the pivot is taken to have minimal $p$-adic valuation,
one can find invertible matrices $N_1 \in \ZZ_q^{r \times r}, N_2 \in \ZZ_q^{s \times s}$ such that
\[ N_1 \cdot A \cdot N_2 \]
is a diagonal matrix whose diagonal elements are called the
\emph{invariant factors} of $A$. We then have the following lemma (the proof is
immediate).

\begin{lemma} \label{invariant}
Let $\theta \in \NN$ be an upper bound for the $p$-adic valuations of the non-zero invariant factors of $A$
and let $N \geq \theta$. Let $\textbf{x}_0 \in \ZZ_q^s$ satisfy
\[ A \cdot \textbf{x}_0 \equiv b \text{ mod } p^N.\]
If there is an $\textbf{x} \in \ZZ_q^s$ such that
\[ A \cdot \textbf{x} = b, \]
then $\textbf{x}$ can be chosen to satisfy $\textbf{x} \equiv \textbf{x}_0$ mod $p^{N-\theta}$.
\end{lemma}

The method works as follows. First, precompute the invariant factors and the matrices $N_1$ and $N_2$ (and their inverses)
modulo $p^{2\theta}$. In total, we need $\widetilde{O}(d^3 n \theta)$ time to do this,
where $d = \max\{r,s\}$ is the dimension of $A$.

Now suppose
we have an $\textbf{x}_0$ such that $A \cdot \textbf{x}_0 \equiv b \text{ mod } p^N$
for some $N \geq \theta$. By Lemma~\ref{invariant}, we can find an $\textbf{x}$ of the form $\textbf{x}_0 + \textbf{t} p^{N-\theta}$
such that $A \cdot \textbf{x} \equiv b \text{ mod } p^{2N}$. To this end, we have to find a $\textbf{t}$ such that
\[ A \cdot \textbf{t} \equiv \frac{b - A \cdot \textbf{x}_0}{p^{N- \theta}} \text{ mod } p^{N + \theta}.\]
Let $T(N)$ denote the time needed to solve a linear system (with fixed linear part $A$) up to precision $N$,
assuming it has a $p$-adic solution. Then
\[ T(2N) = T(N) + T(N + \theta) + \widetilde{O}(d^2nN).\]
Here, the first term comes from the time needed to compute $\textbf{x}_0$.
The last term is dominated by the computation of $A \cdot \textbf{x}_0$ modulo $p^{2N}$.
The second term comes
from the time needed to compute $\textbf{t}$, given $(b - A \cdot \textbf{x}_0)/p^{N - \theta}$ mod $p^{2N}$. Similarly,
$T(N + \theta) = T(N) + T(2\theta) + \widetilde{O}(d^2nN)$. Using our precomputation and the fact that $\theta \leq N$,
we have that $T(2\theta) = \widetilde{O}(d^2nN)$. In conclusion,
\[ T(2N) = 2T(N) + \widetilde{O}(d^2nN).\]
It is obvious that this recurrence relation still holds if $N < \theta$ (again using our precomputation). From a well-known observation in
complexity theory (see for instance \cite[\textsc{Lemma 8.2.}]{Vonzurgathen}) we conclude that
\[ T(N) = \widetilde{O}(d^2nN).\]
Together with our precomputation this results in $\widetilde{O}(d^2nN + d^3n\theta)$ bit-operations.
The following lemma concludes this section.

\begin{lemma} Let $m \in \NN_0$ be the level
at which the reduction starts, i.e. suppose that the polynomial to be reduced is in $S^{(0)}_{[-m,m]}$.
The $p$-adic valuations of the non-zero invariant factors of the matrices $A$ appearing
in our reduction algorithm
are bounded by $\theta = \left\lceil \log_p ( (m + 2(\kappa_2 - \kappa_1 + 1))M + \Delta) \right\rceil$,
where
\[M = \max \{ \pm \ord_P (x) \}_{P \in C \setminus \TT^2_{\Qq}} \quad \text{and} \quad
\Delta = \max \{ - \ord_P (y^{d_t-1}) , - \ord_P (y^{d_b}) \}_{P \in C \setminus \TT^2_{\Qq}}.\]
\end{lemma}

\bproof
We claim that $A$ has the following property: \emph{if $b \in p^{\theta}\ZZ_q^r$
is such that the system $A \cdot \textbf{x} = b$ has a solution in $\QQ_q^s$, then it has a solution in
$\ZZ_q^s$.} Since $N_1$ and $N_2$ are invertible over $\ZZ_q$, this property then still holds for the
matrix $N_1 \cdot A \cdot N_2$, from which the result easily follows.

For simplicity, we will only prove the claim in case
$A$ comes from the system that has to be solved during Step~1 of Phase~1. The other cases work similarly.
Let $b \in p^{\theta}\ZZ_q^r$ be such that $A \cdot \textbf{x} = b$ has a solution in
$\QQ_q^s$. Then $b$ corresponds to a polynomial
\[ h \in S^{(1)}_{[m-c+1, m + 2\kappa_2 - 2]}\]
for which there exists a $g \in S_{[m-c-\kappa_2 + 2, m - 1 + \kappa_2]}$ such that
\[ h - D(g) \in S_{[0,m-c]}.\]
By Corollary~\ref{reductioncoro} and Corollary~\ref{defkappa} (see the first
sentence after the proof of Corollary~\ref{defkappa}), we can reduce this further to eventually obtain a
$g \in S_{[\kappa_1, m - 1 + \kappa_2]}$ such that
\[ h - D(g) \in \mathcal{L}(2D_C).\]
Now, let $\{v_1, \dots, v_m\}$ be a $\QQ_q$-basis for
\[ \frac{\mathcal{L}(2D_C)}{D(\mathcal{L}(D_C))}. \]
As explained in Section~\ref{cohomsec}, this is also a basis for $H^1_{DR}(C)$. In any case,
we can find a $g_0 \in \mL(D_C)$ such that $h - D(g) - D(g_0) = \lambda_1v_1 + \dots + \lambda_mv_m$
for some $\lambda_1, \dots, \lambda_m \in \QQ_q$.

On the other hand, since $h \in S^{(1)}_{[m-c+1, m + 2\kappa_2 - 2]}$ we can
find a $g' \in S^{(0)}_{[\kappa_1, m + 3\kappa_2 - 3]}$ such that
\[ h - D(g') \in \mathcal{L}(2D_C),\]
again by Corollary~\ref{reductioncoro} and Corollary~\ref{defkappa}. Finally, we find a $g_0' \in \mL(D_C)$ for
which $h - D(g') - D(g_0') = \mu_1v_1 + \dots + \mu_mv_m$ for
some $\mu_1, \dots, \mu_m \in \QQ_q$.

Using uniqueness, we conclude that $D(g + g_0) = D(g' + g_0')$. Hence $d(g + g_0) =
\Lambda D(g+g_0) = \Lambda D(g' + g_0') = d(g' + g_0')$ so that
$g + g_0$ and $g' + g_0'$ only differ by a constant.
In particular, $g' \in S^{(0)}_{[\kappa_1, m - 1 + \kappa_2]}$. This concludes the proof.
\eproof


\section{Commode Case}\label{commodesec}

In this section we discuss the simplifications for a
nondegenerate curve with commode Newton polytope.
Note that in practice, this is the most common case.

\begin{definition}
Let $\KK$ be a field. A bivariate polynomial $f \in \KK[\NN^2]$ is called commode if
\[ \forall S \subset \{x, y \} : \dim \Gamma(f_S) = 2 - |S| \, , \]
where $f_S$ denotes the polynomial obtained from $f$ by
setting all variables in $S$ equal to zero.
\end{definition}
The above definition simply means that the Newton polytope $\Gamma(f)$ contains
the origin, a point $(a, 0)$ with $a \in \NN_0$ and a point $(0,b)$ with $b \in \NN_0$.

In the remainder of this section we will assume that $\of \in \Fq[\NN^2]$ is commode and
\nd \ $\Gamma$, in the following sense. A first consequence of the
assumption of commodeness is that $\AB^2_{\Fq}$ is canonically embedded in $\overline{X}_\Gamma$, the toric compactification of
$\TT_{\FF_q}^2$ with respect to $\Gamma$. As such, we can consider $\overline{X}_\Gamma$ as a
compactification of the affine plane, instead of the torus.
Therefore we will work with
a notion of nondegenerateness that is slightly weaker than the one given in Section~\ref{nondegcurves}: it is
no longer necessary to impose the nondegenerateness conditions with respect to the faces lying on the coordinate axes.
However, we now should explicitly impose that $\of$ defines a nonsingular curve in $\AB_{\FF_q}^2$.
For the remainder of this section, we will use this new notion of nondegenerateness. The main geometrical
difference with the old notion
is that now we allow our curve to be tangent to the coordinate axes.
It is also clear that in practice
all elliptic, hyperelliptic and $C_{ab}$ curves can be given by
an equation that is nondegenerate in the above sense.
An important remark is that Corollary~\ref{nullbound}
and Corollary~\ref{liftnondegenerate} still hold under this weaker condition: the proof of Theorem~\ref{Nullstellensatz}
can be adapted to the above situation.

Now, let $\oC$ denote the nonsingular curve $V(\of)$, i.e.\ the closure
in $\overline{X}_\Gamma$ of the locus of $\of$ in $\AB^2_{\Fq}$.
Instead of transforming the curve to the setting described
in Section~\ref{cohomsec}, we will now work with $\of$ itself.
If we furthermore assume that $\of$ is monic in $y$, we obtain similar consequences
as in Section~\ref{cohomsec}, i.e.\
\begin{enumerate}
  \item the set $S := \{ x^k y^l \ | \ k,l \in \NN, \ 0 \leq l < d_y \ \}$ with $d_y = \deg_y \of$ is an
$\FF_q$-basis for $\frac{\FF_q[\NN^2]}{(\of)}$;
  \item every bivariate polynomial in $\Fq[\NN^2]$ has support in $m\Gamma$ for some
big enough $m \in \mathbb{N}$.
\end{enumerate}


\subsubsection{Cohomology of Commode Nondegenerate Curves}

Take an arbitrary lift $f \in \Zq[\NN^2]$ of $\of$ with the same Newton polytope $\Gamma$,
then $f$ is \nd \ $\Gamma$ and $\Gamma$ is commode.  Let $C$ denote the nonsingular
curve obtained by taking the closure of the locus of $f$ in $X_\Gamma$, then
we will compute $H^1_{DR}(C \cap \AB^2_{\Qq})$, instead of $H^1_{DR}(C \cap \TT^2_{\Qq})$.
Note that the difference $C \cap (\AB^2_{\Qq} \setminus \TT^2_{\Qq})$ consists of
$d_x + d_y$ nonsingular points, with $d_x = \deg f(x,0)$ and $d_y = \deg f(0,y)$, which
by Theorem~\ref{comparisontheo} implies that
\[ \dim H^1_{DR}(C \cap \AB^2_{\Qq}) = \dim H^1_{MW}(\oC \cap \AB^2_{\FF_q}) = 2 \Vol(\Gamma) - d_x - d_y + 1 \, . \]

The main difference with the general case is that Theorem~\ref{reduction} needs to be reformulated
as follows, where $A$ is now $\frac{\ZZ_q[\NN^2]}{(f)}$.
\begin{theorem}\label{commodered}
Every element of $D^1(A) \otimes_{\ZZ_q} \QQ_q$
is equivalent modulo $d(A) \otimes_{\ZZ_q} \QQ_q$ with a differential form $\omega$
with divisor
\[ \Div_C(\omega) \geq - D_C - V_C \, , \]
where $V_C$ is defined as $V_C = W_C - (T_x \cap C) - (T_y \cap C)$ with $T_x$ (resp.\
$T_y$) the one-dimensional torus corresponding to the $x$-axis (resp.\ $y$-axis).
\end{theorem}
\bproof
Note that the support of the divisor $D_C$ is disjoint from $(T_x \cap C)$ and from $(T_y \cap C)$
since the corresponding $N_k$ are zero, which explains the definition of $V_C$.  The proof of
Theorem~\ref{reduction} then holds with $W_C$ replaced by $V_C$. Of course, $V_C$ is
still defined over $\QQ_q$.
\eproof

The definition of $\Lambda$ remains the same, but we need to restrict it to
$\mL( - \Div_0(x) - \Div_0(y))$ to obtain a bijection.  Indeed,
\[ \Div_C(\Lambda(h)) = \Div_C(h) + D_C - W_C = \Div_C(h) + D_C - V_C - (T_x \cap C) - (T_y \cap C) \, . \]
Note that $\Div_0(x) = T_x \cap C$ and $\Div_0(y) = T_y \cap C$ and that the support of
$D_C - V_C$ is contained in $C \setminus \AB^2_{\Qq}$; therefore if $\Lambda(h)$
should have no poles on $C \cap \AB^2_{\Qq}$, then clearly $h \in \mL( -\Div_0(x) - \Div_0(y))$.

Theorem~\ref{commodered} implies that each differential form in $D^1(A) \otimes_{\ZZ_q} \QQ_q$
is equivalent modulo exact differential forms with an $\omega$ with divisor
$\Div_C(\omega) \geq - D_C - V_C$.  Let $\omega = \Lambda(h)$, then
\[ \Div_C(\Lambda(h)) \geq -D_C - V_C \Leftrightarrow h \in  \mL(2 D_C - \Div_0(x) - \Div_0(y)) \, . \]
Since the support of $D_C$ is disjoint with the support of $\Div_0(x) + \Div_0(y)$,
we conclude that
\[ \mL(2 D_C - \Div_0(x) - \Div_0(y)) = \mL(2D_C) \cap \mL(- \Div_0(x) - \Div_0(y)) = L_{2 \Gamma}^{-} \, , \]
where $L_{2 \Gamma}^{-}$ denotes the bivariate polynomials with support in $\NN_0^2 \,  \cap \, 2 \Gamma$.
By working modulo $D$ and $f$, we finally obtain the following corollary.

\begin{corollary}
If $\Gamma$ is commode, then $\Lambda$ induces an isomorphism of $\Qq$-vector spaces:
\[ \frac{L^{-}_{2 \Gamma}}{f L^{-}_{\Gamma} + D ( L_\Gamma) } \simeq  H^1_{DR}(C \cap \AB^2_{\Qq}) \, . \]
\end{corollary}


\subsubsection{Lifting Frobenius Endomorphism}  This follows the description given in Section~\ref{Frobeniuslift},
with the simplification that we now only need to compute the action of Frobenius on $x$ and $y$.


\subsubsection{Reduction Algorithm} The reduction in the commode case corresponds to Phase~1 of the
general case as described in Section~\ref{redsec}.  The main difference is that the divisor $E_m$
simplifies to $E_m = m \Div_\infty(x)$; since $f$ is commode, we still have
$S_{[0,m]} \subset \mL(2 D_C + E_m)$, since $E_y = (d_y-1) \Div_\infty(y) \leq D_C$.
Furthermore, if we choose $\kappa$ such that
\[ \mL(D_C + E_y + \Div_\infty(x)) \subset S_{[0,\kappa]} \, , \]
then $\mL(D_C + E_y + E_m) \subset S_{[0, m-1 + \kappa]}$
and $D(S_{[0,m]}) \subset S_{[0,m-1 + \kappa]}$.  The remainder of the algorithm
is then exactly the same with $\kappa_1 = 0$ and $\kappa_2 = \kappa$.


\section{Detailed Algorithm and Complexity Analysis}\label{complsec}

\subsection{Input and output size analysis}
As input our algorithm expects an
$\of \in \FF_q[\ZZ^2]$ ($q=p^n$, $p$ prime) that is \nd \ $\Gamma$, satisfying conditions $1$ and $2$ mentioned
at the beginning of Section~\ref{cohomsec}. A good measure for the input size is
\[ \begin{array}{rcl} \text{number of monomials \ } & \times & \text{ (\ space needed to represent coefficient} \\
                                                 & & \text{ \ + \ space needed to represent exponent vector)}
\\ \end{array}\]
which is $\sim \# (\Gamma \cap \ZZ^2) \cdot (\log q + \log \delta)$, where
$\delta$ is the \emph{degree} of $\of$, that is
\[\max \{ |i|+|j| \, | \, (i,j) \in \Gamma \}.\]
From a result by Scott~\cite{Scott}, that states that
$\#(\Gamma \cap \ZZ^2) \leq 3g + 7$ whenever $g \geq 1$,
it follows that $\# (\Gamma \cap \ZZ^2)$ is asymptotically equivalent
with $g$. Note that the number of points on the boundary $R = \#(\partial \Gamma \cap \ZZ^2)$
is bounded by $2g + 7$.

As output our algorithm gives the characteristic polynomial $\chi(t) := \det (\mathcal{F}^\ast_q - \mathbb{I}t) \in \ZZ[t]$
of the Frobenius morphism $\mathcal{F}^\ast_q$ acting
on $H^1_{MW}(V(\of) \cap \TT^2_{\FF_q} / \QQ_q)$.
A measure for its size follows easily from the Weil conjectures.
Indeed, its degree equals $2 \Vol(\Gamma) + 1$ and $2g$ of its roots have absolute
value $q^{1/2}$. The other roots correspond to $\# (\partial \Gamma \cap \ZZ^2) - 1$ places
lying on $V(f) \setminus \TT^2_{\FF_q}$ and have absolute value $q$. Now, since the $i^\text{th}$ coefficient of
$\chi(t)$ is the sum of $2 \Vol(\Gamma) + 1 \choose i$ $i$-fold products of such roots, we conclude that an
upper bound for the absolute values of the coefficients is given by
\[ {2 \Vol(\Gamma) + 1 \choose \Vol(\Gamma)} q^{g + R - 1} \leq 2^{2 \Vol(\Gamma) + 1} q^{g + R - 1}.\]
Therefore, the number of bits needed to represent $\chi(t)$ is
\[ O\left( (2 \Vol (\Gamma) + 1) \cdot \log (2^{2 \Vol (\Gamma) + 1} q^{g + R - 1}) \right) = O(ng^2)\]
for $p$ fixed.

Note that the zeta function of $V(\of) \cap \TT^2_{\FF_q}$ is then given by
\[ Z_{V(\of) \cap \TT^2_{\FF_q}}(t) = \frac{ \frac{1}{q^{g + R - 1}} \chi(qt)}{1-qt}.\]
The zeta function of the complete model $V(\of)$ can easily be derived from the above.
See \cite{Wouterthesis}
for more details.

\subsection{Asymptotic estimates of some parameters}

We will bound the space and time complexity of our algorithm in terms of $n$ and a set of parameters
that depend only on $\Gamma$ (note that we assume $p$ fixed). In the following,
we will often state that some property holds
for \emph{most common polytopes}: this is not intended to be made mathematically exact. But for instance,
the statement will always hold when $\Gamma$ has a unique right-most and a
unique left-most vertex lying on the $x$-axis, as well as a unique top and a unique bottom vertex lying on the $y$-axis.

The most important parameter is of course $g$,
the number of interior lattice points of $\Gamma$.
During complexity analysis, we can interchange
$g$ with the volume of $\Gamma$ or with the total number of lattice points of $\Gamma$,
as they are all asymptotically equivalent. Indeed,
this follows from Scott's result mentioned above, together with Pick's theorem:
\[ g \leq \#(\Gamma \cap \ZZ^2) \leq 3g + 7 \]
\[ g \leq \text{Vol}(\Gamma) \leq 2g + 3. \]
(given $g \geq 1$). Recall that it follows that $R = (\partial \Gamma \cap \ZZ^2) \leq 2g + 7$. Another parameter is $\delta$, as defined above. We will also make use of the width $w$,
i.e. the maximal difference between the first coordinates of two points of $\Gamma$, and
the height $h$, i.e. $d_t - d_b$. Of course, $h,w \leq 2\delta \leq 2w + 2h$. For most common
polytopes, $wh$ will behave like $g$. However, easy examples show that $w,h$ are in general unbounded for fixed
$g$. For instance, let $\Gamma = \text{Conv}\{ (1,m),(1,m-1),(-1,-m),(-1,-m+1)\}$ for
some arbitrarily big $m \in \NN$. Then $\text{Vol}(\Gamma) = 2$, while $\delta = m + 1$.
See also Remark~\ref{logarithmicterms} below.

Next, we need $\chi_1, \chi_2 \in \ZZ$ such that $\mL(mD_C) \subset S_{[m\chi_1, m\chi_2]}$
for all $m \in \NN_0$. Of course, $\chi_1$ and $\chi_2$ are determined by the
slopes of the top and bottom edges of $\Gamma$. Denote as before the top vertex
with $(c_t, d_t)$ and let $(a,b)$ be the clockwise-next vertex. Suppose that $a \geq c_t$.
Then it is not hard to see that Laurent polynomials with support in the upper
half plane part of $m\Gamma$ reduce (modulo $f$) into $S_{[-\infty, m\tau]}$
where $\tau = c_t + \left\lfloor \frac{d_t(a-c_t)}{d_t - b} \right\rfloor$. Now
\[ \frac{ d_t(a - c_t)}{d_t - b} = (a-c_t) + \frac{ b(a - c_t)}{d_t - b}
\leq w + b(a-c_t) \leq w + 2\text{Vol}(\Gamma) \leq 4g + w + 6.\]
The one but last inequality comes from the fact that the triangle with vertices
$(0,0)$, $(c_t,d_t)$, $(a,b)$ is contained in $\Gamma$. Its volume equals
\[ \frac{ad_t - c_tb}{2} \geq \frac{(a - c_t)b}{2}.\]
Therefore, $\tau \leq c_t + 4g + w + 6$. Using the same argument for the lower half plane,
we conclude that $\mL(mD_C) \subset S_{[-\infty, m(\max(c_t,c_b) + 4g + w + 6)]}$. This
is definitely also true when $a < c_t$. By analogy, $\mL(mD_C) \subset
S_{[m(\min(c_t,c_b) - 4g - w - 6), + \infty]}$, which proves that
we can take $\chi_1, \chi_2$ such that $\chi_2 - \chi_1 \leq 8g + 3w + 12$. For most common polytopes,
$h(\chi_2 - \chi_1)$ is expected to be $O(g^{3/2})$ (by interchanging $x$ and $y$ if necessary).


Strongly related with the foregoing are optimal $\kappa_1, \kappa_2 \in \ZZ$ such
that $\mL(D_C + E_y + \text{Div}_\infty(x) + \text{Div}_0(x)) \subset S_{[\kappa_1, \kappa_2]}$,
see Corollary~\ref{defkappa}.
Note that $\pm \ord_P(x) \leq h$ and $\pm \ord_P(y) \leq w$
for any place $P \in C \setminus \TT^2_{\Qq}$: this follows from Corollary~\ref{inclusion}. Therefore,
$\mL(D_C + E_y + \text{Div}_\infty(x) + \text{Div}_0(x)) \subset \mL((hw + 2h + 1)D_C)$.
By the foregoing, we conclude that we can take $\kappa_2 - \kappa_1 = O(hw(\chi_2 - \chi_1)) = O(hw(g+w))$,
though this is a very rough estimate.
For most common polytopes, a much better bound holds: we can omit $E_y$
(see Remark~\ref{remEm}) and have that $\text{Div}_\infty(x) + \text{Div}_0(x) \leq 2D_C$,
so we can use the same bound as above (multiplied by $3$), i.e. $\kappa_2 - \kappa_1 = O(g + w)$.
Again, for most common polytopes $h(\kappa_2 - \kappa_1)$ is expected to be $O(g^{3/2})$.

Finally, we will often make use of the trivial estimates $g \leq h(\chi_2 - \chi_1), h(\kappa_2 - \kappa_1)$.

\subsection{The algorithm}

\begin{remark} \label{logarithmicterms}
In the introductory section, we mentioned
that the Soft-Oh notation neglects factors
that are logarithmic \emph{in the input size}.
From the example given in the above subsection, it is clear
that factors that are logarithmic in $w,h, \delta, \chi_2 - \chi_1$ and $\kappa_2 - \kappa_1$ need
\emph{not} be logarithmic in the input size.
Nevertheless, we will omit them during complexity analysis. This is mainly for sake of simplicity, but
on the other hand we can prove \cite{Wouterthesis}
that there is some `optimal' setting of $\Gamma$, to be obtained
by stretching and skewing, in which $w, h \sim g^4$. Hence also $\delta, \chi_2 - \chi_1$ and $\kappa_2 - \kappa_1$
are bounded by polynomial expressions in $g$. Moreover, the reduction to this optimal setting goes very fast, as it
is essentially Euclid's algorithm for finding shortest vectors in a lattice.
\end{remark}

\begin{remark} \label{select}
We assume that $\of$ is given as an array of tuples \[(\text{coefficient}, \text{exponent vector})\]
that is ordered with respect to the second components,
so that the coefficient corresponding to a given exponent vector can be selected in $\widetilde{O}(1)$ time.
If this is not the case, this can be easily achieved using a sorting algorithm.
\end{remark}

\noindent \verb"STEP 0: compute" $p$-\verb"adic lift of" $\of$\verb"." First note that we assume
that $\FF_p$ is represented as $\ZZ / (p)$ and that $\FF_q$ is represented as $\FF_p / (\overline{r}(X))$
for some monic irreducible degree $n$ polynomial $\overline{r}(X)$. Take $r(X) \in \ZZ[X]$ such that it
has coefficients in $\{0, \dots, p-1\}$ and reduces to $\overline{r}(X)$ modulo $(p)$. Then
$\ZZ_q$ can be represented as $\ZZ_p / (r(X))$.
Let
\[ \overline{a}_{n-1} [X]^{n-1} + \dots + \overline{a}_{1} [X] + \overline{a}_{0} \]
be any element of $\FF_q$. By the \emph{canonical lift} to $\ZZ_q$, we mean
\[ a_{n-1}[X]^{n-1} + \dots + a_1 [X] + a_0,\]
where the $a_j \in \{0, \dots, p-1\}$ are the unique elements that reduce to $\overline{a}_j$ mod $(p)$.
Finally, if $\of = \sum_{(i,j) \in \ZZ^2 \cap \Gamma} \overline{b}_{ij} x^iy^j$, define
$f = \sum_{(i,j) \in \ZZ^2 \cap \Gamma} b_{ij} x^iy^j$ where the $b_{ij}$ are canonical lifts.\\

\indent \emph{Complexity analysis.} This step needs $\widetilde{O}(ng)$ time and space.\\


\noindent \verb"STEP I: determine" $p$\verb"-adic precision." Assume that all calculations are done
modulo $p^N$ for some $N \in \NN$. What conditions should $N$ satisfy? From the
foregoing, it follows that it suffices to compute $\chi(t)$ modulo $p^{\widetilde{N}}$, where
\[ \widetilde{N} \geq \left\lceil \log_p \left( 2 {2 \Vol(\Gamma) + 1 \choose \Vol(\Gamma)} q^{g + R - 1}  \right) \right\rceil.\]
However, during the reduction process (\verb"STEP V.II") there is some loss of precision: to ensure that everything
remains integral we need to multiply with $p^\varepsilon$ where
\[ \varepsilon = \ceil{ \log_p (m M + \Delta) } \]
with $M = \max \{ \pm \ord_P (x) \}_{P \in C \setminus \TT^2_{\Qq}}$,
$\Delta = \max \{ - \ord_P (y^{d_t-1}) , - \ord_P (y^{d_b}) \}_{P \in C \setminus \TT^2_{\Qq}}$ and
$m = \max \{ |m_1|, |m_2| \}$ the level at which the reduction starts. Here, $m_1, m_2 \in \ZZ$
are such that the objects to be reduced are in $S_{[m_1,m_2]}$.
From Corollary~\ref{inclusion}, it is immediate that $M \leq h$ and $\Delta \leq hw$. To see
what $m$ is bounded by, note that
the objects to be reduced have support in $(9pN + 5p) \Gamma$ (when computed modulo $p^N$).
Indeed, from \verb"STEP V.I" we see that
these objects are of the form
\begin{eqnarray*}                  &   & yf_y \left( \mathcal{F}_p (x^iy^j) \mathcal{F}_p(\beta)
\frac{x \partial \mathcal{F}_p(x)}{\mathcal{F}_p(x)  \partial x} - \mathcal{F}_p(x^iy^j) \mathcal{F}_p(\alpha) \frac{x \partial \mathcal{F}_p(y)}{\mathcal{F}_p(y) \partial x} \right) \\
                                   & - & xf_x \left( \mathcal{F}_p(x^iy^j) \mathcal{F}_p(\beta) \frac{y \partial \mathcal{F}_p(x)}{\mathcal{F}_p(x) \partial y}
- \mathcal{F}_p(x^iy^j) \mathcal{F}_p(\alpha) \frac{y \partial \mathcal{F}_p(y)}{\mathcal{F}_p(y) \partial y} \right),
\end{eqnarray*}
where $(i,j) \in 2\Gamma$. Here $\alpha, \beta \in \ZZ_q[\ZZ^2]$ are Laurent polynomials with support in $2\Gamma$
for which $1 \equiv \alpha x f_x + \beta y f_y$ mod $f$ (see Corollary~\ref{nullbound}). The bound
then follows from Theorem~\ref{thmlift} and the remark below it.

Since $L_{(9pN + 5p)\Gamma} \subset S_{[ (9pN + 5p) \chi_1,
(9pN + 5p) \chi_2]}$, we obtain that
\[ \varepsilon \leq \ceil{ \log_p ((9pN+5p) \max\{ |\chi_1|, \chi_2\} h + hw)}.\]
As a consequence, this is a natural bound on the valuations of the denominators appearing in the matrix of $\mathcal{F}_p^\ast$
(as computed in \verb"STEP VII"). During \verb"STEP VIII" and \verb"STEP IX", our denominators could grow
up to $p^{n(2 \Vol (\Gamma) + 1)\varepsilon}$. In conclusion,
it suffices to take $N$ such that it satisfies $N \geq$
\[ 
  \begin{array}{l} \left\lceil \log_p \left( 2 {2 \Vol(\Gamma) + 1 \choose \Vol(\Gamma)} q^{g + R - 1}  \right) \right\rceil \\
\phantom{} \\
\qquad \qquad \qquad + \ n (2 \Vol(\Gamma) + 1)\ceil{ \log_p ((9pN+5p) \max\{ |\chi_1|, \chi_2\} h + hw) } \\ \end{array} \]
In particular, $N = \widetilde{O}(ng)$.\\


\noindent \verb"STEP II: compute effective Nullstellensatz expansion." In this step, one computes
(up to precision $p^N$) polynomials $\alpha, \beta, \gamma \in \ZZ_q[\ZZ^2]$ with support in $2 \Gamma$ such that
\[ 1 = \gamma f + \alpha x \frac{\partial f}{\partial x} + \beta y \frac{\partial f}{\partial y}.\]
This defines a linear system $A \cdot \textbf{x} = B$ that can be solved using Gaussian elimination, in each step of
which the pivot is taken to be a $p$-adic unit. This is possible since
the linear map defined by $A$ is surjective (by Theorem~\ref{Nullstellensatz}). In particular,
there is no loss of precision. Note that instead of
Gaussian elimination, one can use the method described at the end of Section~\ref{redsec}. In this way,
one gains a factor $g$ time. But for the overall complexity analysis this makes no difference.\\


\indent \emph{Complexity analysis.}
Selecting the entries of $A$ takes $\widetilde{O}(g^2)$ time (see Remark~\ref{select}). One then
needs $\widetilde{O}(nNg^3) = \widetilde{O}(n^2g^4)$
time and $O(nNg^2) = \widetilde{O}(n^2g^3)$ space to solve the system.\\

\noindent \verb"STEP III: compute lift of Frobenius." Take lifts $\delta, \delta_x, \delta_y \in \ZZ_q[\ZZ^2]$ of $\overline{\gamma}^p$,
$\overline{\alpha}^p$, $\overline{\beta}^p$ and compute a zero of the polynomial
\[ H(Z) = (1 + \delta_xZ)^a(1 + \delta_yZ)^bf^\sigma (x^p(1 + \delta_xZ), y^p(1 + \delta_yZ))\]
(as described in Section~\ref{Frobeniuslift}) up to precision $p^N$, using Newton iteration and
starting from the approximate solution 0.
Reduce all intermediate calculations
modulo $f$ to the basis $\{ \, x^iy^j \, | \, d_b \leq j < d_t \}$ (this is why
the terms $- (a \delta_x + b \delta_y - \delta)f^pZ - f^p$, that were added for theoretical reasons,
can be omitted in the formula for $H(Z)$).
Finally, if we denote the result by $Z_0$, expand $Z_x := 1 + \delta_xZ_0$, $Z_y := 1 + \delta_yZ_0$
and compute their inverses up to precision $p^N$ using Newton iteration (again reduce the intermediate
calculations modulo $f$). Note that if we take $a$ and $b$ minimal, then $\deg H \leq w + h$.\\

\indent \emph{Complexity analysis.}
Remark that it is better \emph{not} to expand the polynomial $H(Z)$ (nor its derivative
$\frac{dH}{dZ}(Z)$), but to leave it in the above compact representation. The reason is that
the expanded versions of $H$ and $\frac{dH}{dZ}$ are very space-costly.

A similar complexity estimate has been made in \cite{DVCab}. The complexity
is dominated by the last iteration step, which in its turn is dominated by $O(g)$ computations of
terms of the form
\[ (1 + \delta_xZ')^i(1 + \delta_yZ')^j \]
where $Z' \in S_{[6pN\chi_1, 6pN\chi_2]}$, $i \in \{0, \dots, w\}$ and $j \in \{0, \dots, h\}$ (because of
(\ref{boundforrootsofH})).
Note that reducing a polynomial with support in $[6pN\chi_1, 6pN\chi_2]\times[-\lambda d_b, \lambda (d_t-1)]$
(for some $\lambda \in \NN_0$) to the basis mentioned above can be done in $\widetilde{O}(\lambda h N (\chi_2 - \chi_1) \cdot g \cdot nN) =
\widetilde{O}(\lambda n^3g^3h(\chi_2 - \chi_1))$
time (at least if we know that
all intermediate results are supported in $[6pN\chi_1, 6pN\chi_2] \times \ZZ$ modulo $p^N$). Therefore,
the overall time complexity of \verb"STEP III" amounts to $\widetilde{O}(n^3g^4h(\chi_2 - \chi_1))$, whereas the
space complexity is $\widetilde{O}(n^3g^2h(\chi_2 - \chi_1))$. Note that this indeed dominates the
time and space needed to compute the Frobenius substitutions, each of which can be done in $\widetilde{O}(n \cdot nN)$
time (see e.g. \cite[Section 12.5]{Handbook}).

The complexity of computing $Z_x, Z_y, Z_x^{-1}, Z_y^{-1}$ works similarly and is dominated by the above.\\

\noindent \verb"STEP IV: `precompute'" $\mathcal{F}_p^\ast(dx/xyf_y)$\verb"." Here, $\mathcal{F}_p^\ast$ is the
$\QQ_q$-vector space endomorphism of $\Omega_C(C \cap \TT^2_{\Qq})$ induced by $\mathcal{F}_p$.
Note that $dx / f_y = \beta y dx - \alpha x dy$.  Thus $\mathcal{F}_p^\ast(dx/xyf_y) =$
\[ \mathcal{F}_p(\beta) \left( \frac{\partial \mathcal{F}_p(x)}{\mathcal{F}_p(x) \partial x}dx
+ \frac{\partial  \mathcal{F}_p(x)}{\mathcal{F}_p(x) \partial y}dy \right) - \mathcal{F}_p(\alpha)
\left( \frac{\partial \mathcal{F}_p(y)}{\mathcal{F}_p(y) \partial x}dx
+ \frac{\partial  \mathcal{F}_p(y)}{\mathcal{F}_p(y) \partial y}dy \right)\]
\[ = \left( \mathcal{F}_p(\beta) \frac{ \partial \mathcal{F}_p(x)}{\mathcal{F}_p(x)\partial x} -
\mathcal{F}_p(\alpha) \frac{ \partial \mathcal{F}_p(y)}{\mathcal{F}_p(y)\partial x} \right)dx
+ \left( \mathcal{F}_p(\beta) \frac{ \partial \mathcal{F}_p(x)}{\mathcal{F}_p(x)\partial y} -
\mathcal{F}_p(\alpha) \frac{ \partial \mathcal{F}_p(y)}{\mathcal{F}_p(y)\partial y} \right) dy \]
However, as will become clear in the following step, it is more natural to precompute
\[ E := yf_y \left( \mathcal{F}_p(\beta) \frac{ x \partial \mathcal{F}_p(x)}{\mathcal{F}_p(x)\partial x} -
\mathcal{F}_p(\alpha) \frac{ x \partial \mathcal{F}_p(y)}{\mathcal{F}_p(y)\partial x} \right) - xf_x \left( \mathcal{F}_p(\beta) \frac{ y \partial \mathcal{F}_p(x)}{\mathcal{F}_p(x)\partial y} -
\mathcal{F}_p(\alpha) \frac{ y \partial \mathcal{F}_p(y)}{\mathcal{F}_p(y)\partial y} \right). \]
Furthermore, this object has nicer convergence properties, in the
sense that it is supported modulo $p^N$ in an easy to determine multiple of $\Gamma$ ($(9pN + 3p)\Gamma$ to
be precise). Therefore, we have a good control (in terms of $\chi_1$ and $\chi_2$) on the
size of the objects we are computing with.\\

\indent \emph{Complexity analysis.} The complexity of this step is dominated by
the computation of $O(g)$ expressions of the form $Z_x^iZ_y^j$, where $|i|$ and $|j|$
are $O(\delta)$. As before, this results in $\widetilde{O}(n^3g^4h(\chi_2 - \chi_1))$ time and
$\widetilde{O}(n^3g^2h(\chi_2 - \chi_1))$ space.\\

\noindent \verb"STEP V: for every" $(i,j) \in 2 \Gamma$\verb":"\\

\noindent \verb"STEP V.I: let Frobenius act on" $x^iy^j$\verb"." In this step, one actually computes
\[\Lambda^{-1}(\mathcal{F}_p^\ast(\Lambda(x^iy^j))).\]
Note that $\mathcal{F}_p^\ast(\Lambda(x^iy^j))$ is given by
$\mathcal{F}_p(x^iy^j) \mathcal{F}_p^\ast(dx/xyf_y)$.
To translate back, if
\[ \mathcal{F}_p^\ast( \Lambda(x^iy^j)) = g_{ij,1} dx + g_{ij,2} dy\]
then
\[ \Lambda^{-1}(\mathcal{F}_p^\ast(\Lambda(x^iy^j))) = xy(f_y g_{ij,1} - f_x g_{ij,2}).\]
Therefore, we output
\[ \mathcal{F}_p(x^iy^j) \cdot E \]
where $E$ is the expression that was precomputed during the foregoing step.\\

\noindent \verb"STEP V.II: reduce modulo" $D.$ Using the method described in Section~\ref{redsec}, reduce
the output of the foregoing substep (after multiplying with $p^\varepsilon$) to obtain polynomials
$r_{ij} \in \mathcal{L}^{(1)}(2D_C) \subset L^{(0)}_{2\Gamma}$. Note that we want our
output $r_{ij}$ to be supported in $2\Gamma$: at this stage, we are no longer interested in the reduction to the basis
$\{ \, x^iy^j \, | \, d_b \leq j < d_t \}$.\\

\indent \emph{Complexity analysis.} The complexity of the first substep can
be estimated using a method similar to what we did in \verb"STEP IV",
resulting in $\widetilde{O}(n^3g^3h(\chi_2 - \chi_1))$ time (per monomial) and
$\widetilde{O}(n^3g^2h(\chi_2 - \chi_1))$ space. For the second substep, it suffices to
analyze the complexity of Phase~1 and Step~{t+1}, as described in Section~\ref{redsec}.
During Phase~1, one needs to solve systems of size $\sim h(2\kappa_2 + c)$. Therefore, it
is optimal to choose $c = \kappa_2$. The number of systems to be solved is then
bounded by $m/c = m/\kappa_2$. Using similar estimates for Phase 2 and using the analysis made at
the end of Section~\ref{redsec},
this results in a use of
\[ \widetilde{O}(h^2 (\kappa_2 - \kappa_1)(\chi_2 - \chi_1) n N^2 + h^3(\kappa_2 - \kappa_1)^2(\chi_2 - \chi_1)nN) \]
time before proceeding to Step~{t+1}. In this final step, one needs to solve a linear system of size
$O(h \max \{\kappa_2 - \kappa_1, \chi_2 - \chi_1\})$, resulting in a time-cost of
\[ \widetilde{O}(h^2 \left( \max\{ \kappa_2 - \kappa_1, \chi_2 - \chi_1 \} \right)^2 n N +
h^3 \left( \max\{ \kappa_2 - \kappa_1, \chi_2 - \chi_1 \} \right)^3 n).\]
The extra space needed during Phase~1 and Step~{t+1} is
\[\widetilde{O}(h^2 \left( \max\{ \kappa_2 - \kappa_1, \chi_2 - \chi_1 \} \right)^2 n N),\]
though this will in general be dominated by the space needed to store the polynomial $h$ that is to be reduced,
which is $\widetilde{O}(n^3g^2h(\chi_2 - \chi_1))$.

Since substeps \verb"V.I" and \verb"V.II" have to be executed for $O(g)$ monomials,
we obtain the following global estimates for \verb"STEP V": a time-cost of
\[ \widetilde{O}(n^3g^3h^2 \left( \max\{ \kappa_2 - \kappa_1, \chi_2 - \chi_1 \} \right)^2
    + n^2g^2 h^3 \left( \max\{ \kappa_2 - \kappa_1, \chi_2 - \chi_1 \} \right)^3) \]
and a space-cost of $\widetilde{O}(n^3gh^2\left( \max\{ \kappa_2 - \kappa_1, \chi_2 - \chi_1 \} \right)^2)$.

Note that our time-estimate dominates the time needed to actually \emph{compose} the systems that
are to be solved.\\



\noindent \verb"STEP VI: compute a" $\ZZ_q$\verb"-basis of" $M_H = \frac{\mathcal{L}^{(0)}(2D_C)}{\left(D(\mathcal{L}(D_C))\right)^{(0)}}$.
Note that from the proof of Lemma~\ref{strongRRSpace}, we have that $\mathcal{L}^{(0)}(mD_C) = L^{(0)}_{m\Gamma}$
for any $m \in \NN_0$. Therefore, we actually have to compute a $\ZZ_q$-basis of
\[ \frac{L^{(0)}_{2\Gamma}}{\left( D(L_\Gamma) + f L_\Gamma \right)^{(0)}}.\]
Consider the module $D(L^{(0)}_\Gamma) + f L^{(0)}_\Gamma$
and express a vector $\mathcal{A}$ whose entries are the generators
$\left\{  D(x^iy^j), fx^iy^j \right\}_{(i,j) \in \Gamma \cap \ZZ^2}$
in terms of a vector $\mathcal{B}$ whose entries are $\left\{ x^ry^s \right\}_{(r,s) \in 2\Gamma \cap \ZZ^2}$:
\[ \mathcal{A} = E \cdot \mathcal{B}.\]
Now compute $\ZZ_q$-invertible matrices $N_1$ and $N_2$ (and their inverses) such that
$N_1 \cdot E \cdot N_2$ is a diagonal matrix. Its non-zero entries are the non-zero
invariant factors of $E$ and will be denoted by $d_1, \dots, d_\ell$.
If we write
\[ N_1 \cdot \mathcal{A} = N_1 \cdot E \cdot N_2 \cdot N_2^{-1} \cdot \mathcal{B}, \]
we see that the entries of $N_2^{-1} \cdot \mathcal{B}$ form a basis
$\{f_1, \dots, f_k\}$ of $L^{(0)}_{2\Gamma}$ such that $\{d_1f_1, \dots, d_\ell f_\ell\}$
is a basis of $D(L^{(0)}_\Gamma) + f L^{(0)}_\Gamma$. It is then easily seen that
$\{f_1, \dots, f_\ell\}$ is a basis of $\left( D(L_\Gamma) + f L_\Gamma \right)^{(0)}$.
Finally, $\{f_{\ell + 1}, \dots, f_k\}$ is a basis of $M_H$.

When computing modulo a finite precision, some caution is needed: to determine $f_{\ell + 1}, \dots, f_k$ modulo $p^N$,
it does not suffice to do the above computations modulo the same precision. During this step (and only during this step),
we need to compute modulo $p^{N + N_0}$, where $N_0 = \lfloor \ell n \log_p(\ell whnp) \rfloor + 1 = O(N)$. Indeed, we
claim that $N_0$ is a strict upper bound for the $p$-adic valuation of any non-zero $(\ell \times \ell)$-minor of $E$. As a consequence,
the valuations of the non-zero invariant factors $d_1, \dots, d_\ell$ are also strictly bounded by $N_0$.
Therefore, we will be
able to find invertible matrices $\widetilde{N}_1$ and $\widetilde{N}_2$ such that
\[ \widetilde{N}_1 \cdot E \cdot \widetilde{N}_2^{-1} \]
is congruent modulo $p^{N + N_0}$ to the above diagonal matrix. The `basis' $\{ \widetilde{f}_{\ell + 1}, \dots, \widetilde{f}_k \}$
we find in this way
corresponds modulo $p^N$ to the basis mentioned above: if we would want to finalize the above
diagonalization (which was only carried out modulo $p^{N+N_0}$), we would need to subtract from the $\widetilde{f}_i$
Laurent polynomials with coefficients divisible by $p^{(N + N_0)}/p^{N_0} = p^N$. Actually, one can check that
$\{ \widetilde{f}_{\ell + 1}, \dots, \widetilde{f}_k \}$ is a basis itself, but we won't need this.
If in \verb"STEP VII" we write $f_{\ell + 1}, \dots, f_k$ and $N_2$, we actually mean
the reductions mod $p^N$ of $\widetilde{f}_{\ell + 1}, \dots, \widetilde{f}_k$ and $\widetilde{N}_2$ that
were computed this way.

It remains to prove the claim, i.e. the $p$-adic valuation of any non-zero $(\ell \times \ell)$-minor of $E$
is strictly bounded by $N_0$. Let $r(X)$ be the polynomial from \verb"STEP 0" and let $\theta \in \CC$ be
a root of it. Consider $K = \QQ(\theta)$ and let $\mathcal{O}_K$ be its ring of algebraic integers.
Then $\mathfrak{p} = (p) \subset \mathcal{O}_K$ is a prime ideal and the $\mathfrak{p}$-adic completion
of $K$ can be identified with $\QQ_q$. Under this identification, $E$ has entries
\[ \sum_{i=0}^{n-1} a_i \theta^i \in \mathcal{O}_K\]
where the $a_i \in \ZZ$ satisfy $|a_i| \leq 2whp$. Since the complex norm of
\emph{any} root of $r(X)$ is bounded by $p$ by Cauchy's bound, we conclude that the entries $e$ of $E$
satisfy
\[ |e_{ij}|_K \leq nwhp^n \leq (whnp)^n\]
for \emph{any} archimedean norm $| \cdot |_K$ on $K$ that extends the classical absolute value on $\QQ$. Since an $(\ell \times \ell)$-minor $m$
is the sum of $\ell!$ $\ell$-fold products of such entries, it follows that
\[ |m|_K \leq (\ell whnp)^{\ell n}.\]
Since $m$ is an algebraic integer, from the product formula we have
\[ | m |_{\mathfrak{p}}^{-n} \leq \prod | m|_K \leq (\ell whnp)^{\ell n^2} \]
(if $m \neq 0$), where $| \cdot |_{\mathfrak{p}}$ is scaled such that $|p|_{\mathfrak{p}} = 1/p$ and where the product is over all archimedean norms $|\cdot|_K$ on $K$, to be counted
twice if it comes from a non-real root of $r(X)$. From this we finally get that $\ord_p m \leq \ell n \log_p(\ell whnp)$.\\


\indent \emph{Complexity analysis.}
This step needs $O(g^3)$ ring operations, each of which takes
$\widetilde{O}(nN)$ time. Therefore, the time complexity of this step is $\widetilde{O}(n^2g^4)$
while the space complexity amounts to $\widetilde{O}(n^2g^3)$.\\

\noindent \verb"STEP VII: compute a matrix of" $p$\verb"-th power Frobenius."
From \verb"STEP V", we know that $p^\varepsilon x^iy^j$ is mapped to $r_{ij}$. Therefore, it is
straightforward to compute the action of Frobenius on $f_{\ell + 1}, \dots, f_k$ and express
it in terms of $\mathcal{B}$:
\[ \Lambda^{-1} \mathcal{F}^\ast_p \Lambda \ p^\varepsilon  \left( \begin{array}{c} f_{\ell + 1} \\ \vdots \\ f_k \\ \end{array} \right) = F \cdot \mathcal{B}.\]
Since $F \cdot \mathcal{B} = F \cdot N_2 \cdot N_2^{-1} \cdot \mathcal{B}$, we obtain a matrix of Frobenius
as $p^{-\varepsilon}$ times an appropriate submatrix $M$ of $F \cdot N_2$.\\

\indent \emph{Complexity analysis.}
The complexity of this step is dominated by the computation of $F \cdot N_2$,
which takes $\widetilde{O}(n^2g^4)$ time and $\widetilde{O}(n^2g^3)$ space,
and by $O(g^2)$ Frobenius substitutions, taking an extra
$\widetilde{O}(g^2 \cdot n \cdot nN) = \widetilde{O}(n^3g^3)$ time.\\

\noindent \verb"STEP VIII: compute a matrix of" $q$\verb"-th power Frobenius."
The matrix $p^{-\varepsilon}M$ of the foregoing step
is a matrix of $\mathcal{F}_p^\ast$, which is a $\QQ_p$-vector
space morphism acting on $H^1_{MW}(C \cap \TT_{\QQ_q}^2)$. A
matrix of $\mathcal{F}_q^\ast$ is then given by
$p^{-n\varepsilon}\mathcal{M}_n$ where $\mathcal{M}_n =
M^{\sigma^{n-1}}\cdot M^{\sigma^{n-2}} \cdots M^\sigma \cdot M$.


$\mathcal{M}_n$ can be computed using the following method that was presented by Kedlaya \cite{Kedlaya}: let $n =
\mathfrak{n}_1\mathfrak{n}_2\dots \mathfrak{n}_k$
be the binary expansion of $n$ and write $n' = \mathfrak{n}_1\mathfrak{n}_2\cdots\mathfrak{n}_{k-1}$, then we
have the formula
\[ \mathcal{M}_n = \mathcal{M}_{n'}^{\sigma^{n' + \mathfrak{n}_k}} \cdot
\mathcal{M}_{n'}^{\sigma^{\mathfrak{n}_k}} \cdot M^{\mathfrak{n}_k}\]
by means of which $\mathcal{M}_n$ can be computed recursively .\\

\indent \emph{Complexity analysis.} Applying some $\sigma^i$ ($i \leq n$) to a matrix of size $O(g)$ takes
$\widetilde{O}(g^2\cdot n \cdot nN) = \widetilde{O}(n^3g^3)$ time, if we
precompute $[X]^{\sigma^i}$ as a root of the polynomial $r$ that defines $\ZZ_q$, using
Newton iteration and starting from the approximate solution $[X]^{p^i} \in \FF_q$.
The complexity of \verb"STEP VIII" is then dominated by $O(\log n)$ matrix multiplications
and $O(\log n)$ applications of some $\sigma^i$, resulting in $\widetilde{O}((n+g)n^2g^3)$ time.
The space needed is $\widetilde{O}(n^2g^3)$. \\

\noindent \verb"STEP IX: output the characteristic polynomial of Frobenius".
The characteristic polynomial $\widetilde{\chi}(t)$ of $\mathcal{M}_n$
can be computed using the classical algorithm based on the
reduction to the Hessenberg form \cite[Section 2.4.4]{Cohen}. In each step of this reduction,
the pivot should be chosen to be an entry under the diagonal with minimal
$p$-adic valuation (unless this exceeds the required precision). In this way, no denominators are introduced.
Write
\[ \widetilde{\chi}(t) = \sum_{i=0}^{2\Vol(\Gamma) + 1} c_it^i.\]
Then the characteristic polynomial of $\mathcal{F}_q^\ast$ (or of
$p^{-n\varepsilon}\mathcal{M}_n$) is given by
\[ \chi(t) = \sum_{i=0}^{2\Vol(\Gamma) + 1} p^{(i-2\Vol(\Gamma) - 1)n\varepsilon}c_it^i \quad \in \ZZ[t].\]
This finalizes the description of the algorithm.\\

\indent \emph{Complexity analysis.} This needs $O(g^3 \cdot nN) = \widetilde{O}(n^2g^4)$
time and $\widetilde{O}(n^2g^3)$ space.

\subsection{Main theorem}

The above analysis allows us to reformulate Theorem~\ref{main} in more detail.

\begin{theorem} \label{maintheorem}
There exists a deterministic algorithm to compute the zeta function of a
bivariate Laurent polynomial $\overline{f} \in \mathbb{F}_{p^n}[\ZZ^2]$ that
is nondegenerate with respect to its Newton polytope $\Gamma$, given that the latter
contains the origin and has unique top and bottom vertices. Let $g, h, w, \kappa_1, \kappa_2, \chi_1, \chi_2$
be as above. Then for fixed $p$, it has running time
\[ \widetilde{O}(n^3g^3h^2 \left( \max\{ \kappa_2 - \kappa_1, \chi_2 - \chi_1 \} \right)^2
    + n^2g^2 h^3 \left( \max\{ \kappa_2 - \kappa_1, \chi_2 - \chi_1 \} \right)^3). \]
The space
complexity amounts to
\[ \widetilde{O}(n^3gh^2\left( \max\{ \kappa_2 - \kappa_1, \chi_2 - \chi_1 \} \right)^2)\]
The $\widetilde{O}$-notation hides factors that are logarithmic in $n,g,w$ and $h$.
For `most common' polytopes, the estimates
$h(\chi_2 - \chi_1) \approx h(\kappa_2 - \kappa_1) \approx g^{3/2}$
hold, so that the algorithm needs $\widetilde{O}(n^3g^6 + n^2g^{6.5})$ time and $\widetilde{O}(n^3g^4)$
space.
\end{theorem}

Recall from Section~\ref{cohomsec} that the above conditions on $\Gamma$ are not restrictive.
Note that in the $C_{ab}$ curve case, a better estimate
for $h(\chi_2 - \chi_1) = h(\kappa_2 - \kappa_1)$ is $g$, yielding
a time complexity of $\widetilde{O}(n^3g^5)$ and a space
complexity of $\widetilde{O}(n^3g^3)$. This is the same as in the algorithm presented in \cite{DVCab}.


\section{Conclusions}\label{conclsec}

In this paper, we presented a generalization of Kedlaya's algorithm
to compute the zeta function of a nondegenerate curve over a finite field of small characteristic.
As the condition of nondegenerateness is generic, the algorithm works for curves that
are defined by a randomly chosen bivariate Laurent polynomial with given Newton polytope $\Gamma$.
It requires $\widetilde{O}(n^3 \Psi_t)$ amount of time and $\widetilde{O}(n^3 \Psi_s)$ amount
of space, where $\Psi_t, \Psi_s$ are functions that depend on $\Gamma$ only. For non-exotic choices of
$\Gamma$, we have that $\Psi_t \sim g^{6.5}$ and $\Psi_s \sim g^{4}$, where $g$ is the number of interior
lattice points of $\Gamma$ (which is precisely the geometric genus of the curve).
(in fact, if $n \gg g$, which will usually be the case, we have $\Psi_t \sim g^6$).
In the case of a $C_{ab}$ curve, we obtain the estimates $\Psi_t \sim g^5$ and
$\Psi_s \sim g^3$, so that the algorithm works (at least asymptotically) as fast as
the one presented in \cite{DVCab}.
At this moment, the algorithm has not yet been fully implemented.

In order to develop the algorithm, we proved a number of theoretical results on nondegenerate curves that
are interesting in their own right, for example a linear effective Nullstellensatz
for sparse Laurent polynomials in any number of variables. Also, we adapted the Frobenius lifting technique used
in \cite{DVCab} to prove a convergence rate in which the Newton polytope $\Gamma$ plays a very natural role.

These results seem to reveal an entirely sparse description of the first
Monsky-Washnitzer cohomology group and the action of Frobenius on it,
though this should be investigated further.
In particular, during reduction modulo exact differentials
we loose track of the Newton polytope for complexity reasons.


\end{document}